\documentclass[12pt]{article}
\usepackage{color,latexsym,amsfonts,amssymb}

\usepackage{amssymb,amsmath}

\begin{document}

\baselineskip15pt

\font\fett=cmmib10
\font\fetts=cmmib7
\font\bigbf=cmbx10 scaled \magstep2
\font\bigrm=cmr10 scaled \magstep2
\font\bbigrm=cmr10 scaled \magstep3
\def\cnl{\centerline}
\def\etal{{\it et al.}}
\def\text{{}}
\def\un{^{(n)}}
\def\ff{{\cal F}}
\def\Ref#1{(\ref{#1})}

\def\Blm{\left|}
\def\Brm{\right|}
\def\Bl{\left(}
\def\Br{\right)}
\def\nti{n\to\infty}
\def\lnti{\lim_{\nti}}
\def\law{{\cal L}}
\def\sji{\sum_{j\ge1}}
\def\Cal{\cal}

\def\real{\text{\rm I\kern-2pt R}}
\def\re{\real}
\def\expec{\text{\rm I\kern-2pt E}}
\def\ex{\expec}
\def\prob{\text{\rm I\kern-2pt P}}
\def\pr{\prob}
\def\rat{\text{\rm Q\kern-5.5pt\vrule height7pt depth-1pt\kern4.5pt}}
\def\comp{\text{\rm C\kern-4.7pt\vrule height7pt depth-1pt\kern4.5pt}}
\def\nat{\text{\rm I\kern-2pt N}}
\def\integ{{\bf Z}}
\def\qedbox{\vcenter{\hrule height.5pt\hbox{\vrule width.5pt height8pt
\kern8pt\vrule width.5pt}\hrule height.5pt}}
\def\half{{\textstyle {1 \over 2}}}
\def\quarter{{\textstyle {1 \over 4}}}
%\mathchardef\scrl="024C
%\mathchardef\scrd="0244
\def\scrl{{\Cal L}}
\def\scrd{{\Cal D}}
\def\tod{\buildrel \scrd \over \longrightarrow}
\def\eqd{\buildrel \scrd \over =}
\def\tolone{\buildrel L_1 \over \longrightarrow}
\def\toltwo{\buildrel L_2 \over \longrightarrow}
\def\tolp{\buildrel L_p \over \longrightarrow}
\def\l{\lambda}
\def\L{\Lambda}
\def\a{\alpha}
\def\b{\beta}
\def\g{\gamma}
\def\G{\Gamma}
\def\d{\delta}
\def\D{\Delta}
\def\e{\varepsilon}
\def\h{\eta}
\def\z{\zeta}
\def\th{\theta}
\def\k{\kappa}
\def\m{\mu}
\def\n{\nu}
\def\p{\pi}
\def\r{\rho}
\def\s{\sigma}
\def\S{\Sigma}
\def\t{\tau}
\def\f{\varphi}
\def\ch{\chi}
\def\ps{\psi}
\def\o{\omega}
\def\lee{\,\le\,}
\def\leee{\quad\le\quad}
\def\gee{\,\ge\,}
\def\geee{\quad\ge\quad}
\def\scrn{{\Cal N}}
\def\scra{{\Cal A}}
\def\scrf{{\Cal F}}
\def\var{\text{\rm Var\,}}
\def\cov{\text{\rm Cov\,}}

\def\sin{\sum_{i=1}^n}
\def\sn{\sum_{i=1}^n}
\def\cross{\times}

\def\nin{\noindent}
\def\bsk{\bigskip}
\def\msk{\medskip}
\def\widebar{\bar}

\def\sno{\sum_{n\ge 0}}
\def\sj{\sum_{j\ge 0}}
\def\proof{\noindent{\bf Proof.}\ }
\def\remark #1. {\medbreak\noindent{\bf Remark #1.}\quad}
\def\remarks{\medbreak\noindent{\bf Remarks.}\quad}
\def\head #1/ {\noindent{\bf #1.}\medskip}
\def\dtv{d_{TV}}
\def\Po{\text{\rm Po\,}}
\def\Bi{\text{\rm Bi\,}}
\def\Be{\text{\rm Be\,}}
\def\CP{\text{\rm CP\,}}

\def\pb{\hbox{{\fett\char'031}}}
\def\pbh{\text{\bf{\hbox{\^{\fett\char'031}}}}}
\def\bp{\pb}
\def\lb{\hbox{{\fett\char'025}}}
\def\lbs{\hbox{{\fetts\char'025}}}
\def\bl{\lb}
\def\rb{\hbox{\fett\char'032}}
\def\sb{\hbox{\fett\char'033}}
\def\iid{\text{\rm independent\ and\ identically\ distributed}}

\def\and{\ \text{\rm{ and }}\ }
\def\for{\ \text{\rm{ for }}\ }
\def\forall{\ \text{\rm{ for all }}\ }
\def\forsome{\ \text{\rm{ for some }}\ }
\def\tif{\text{{\rm if\ }}}
\def\tin{\text{{\rm in\ }}}
\def\lcm{\text{\rm{l.c.m.\,}}}

\def\xx{{\cal X}}
\def\supp{{\rm supp\,}}

\def\BCL{Barbour, Chen \& Loh}
\def\ABT{Arratia, Barbour \& Tavar\'e}
\def\BHJ{Barbour, Holst \& Janson}

\def\Blb{\left\{}
\def\Brb{\right\}}

\def\giv{\,|\,}
\def\Giv{\,\Big|\,}
\def\ep{\hfill$\Box$}
\def\bone{{\bf 1}}
\def\non{\nonumber}

\def\adb{}
\def\sodmi{\sum_{s=0}^{d-1}}
\def\sie{\sum_{i\in E}}
\def\sii{\sum_{i\ge1}}
\def\Le{\ \le\ }
\def\TP{{\rm TP\,}}
\def\Ord#1{O\Bl #1 \Br}
\newcommand{\eqs}{\begin{eqnarray*}}
\newcommand{\ens}{\end{eqnarray*}}

\newcommand{\eqa}{\begin{eqnarray}}
\newcommand{\ena}{\end{eqnarray}}
\newcommand{\eq}{\begin{equation}}
\newcommand{\en}{\end{equation}}

\def\numberlikeadb{\global\def\theequation{\thesection.\arabic{equation}}}
\numberlikeadb
\newtheorem{theorem}{Theorem}[section]
\newtheorem{lemma}[theorem]{Lemma}
\newtheorem{corollary}[theorem]{Corollary}
\newtheorem{proposition}[theorem]{Proposition}
\newtheorem{example}[theorem]{Example}

\def\pr{{\mathbb P}}
\def\ex{{\mathbb E}}
\def\re{{\mathbb R}}
\def\integ{{\mathbb Z}}
\def\ZZ{\integ}
\def\Ceka{\v Cekan\-avi\v cius}
\def\sjom{\sum_{j=0}^m}
\def\uo{^{(0)}}
\def\ul{^{(l)}}
\def\gg{{\mathcal G}}
\def\bh{{\bar h}}
\def\bX{{\widebar X}}
\def\bU{{\widebar U}}
\def\BX{Barbour \&~Xia}
\def\BC{Barbour \&  \Ceka}
\def\sro{\sum_{r\ge0}}
\def\srz{\sum_{r\in\integ}}
\def\ssz{\sum_{s\in\integ}}
\def\diff{\lfloor \m-\s^2 \rfloor}
\def\diffi{\lfloor \m_1-\s_1^2 \rfloor}
\def\difft{\lfloor \m_2-\s_2^2 \rfloor}
\def\sid{\s_1^2 + \d_1}
\def\std{\s_2^2 + \d_2}
\def\hX{{\widehat X}}
\def\hY{{\widehat Y}}
\def\hW{{\widehat W}}
\def\hQ{{\widehat Q}}
\def\hS{{\widehat S}}
\def\Ge{{\rm Ge\,}}
\def\sjn{\sum_{j=1}^n}
\def\bone{{\bf 1}}
\def\ixo{\bone_{\xx_0}}
\def\ixoc{{\bf 1}_{\xx_0^c}}
\def\aaomi{\aaa_0^{-1}}
\def\uu{{\cal U}}
\def\Eq{\ =\ }
\def\bb{{\cal B}}
\def\aaa{{\cal A}}
\def\sjo{\sum_{j\ge0}}
\def\dhh{d_{\hh}}
\def\dff{d_{\bff}}
\def\ngg{\|_{\gg}}
\def\bff{{\overline\ff}}
\def\nsup{\|_\infty}
\def\thalf{\tfrac12}
\def\derivbnd{\sqrt{2\pi}}
\def\sli{\sum_{l\ge1}}
\def\nn{{\cal N}}
\def\hh{{\cal H}}
\def\ui{^{(1)}}
\def\ignore#1{}
\def\xxo{{\bf X}}
\def\LL{{\overline L}}
\def\tl{{\tilde\l}}
\def\Ge{\ \ge\ }
\def\Lip{{\rm Lip}}
\def\dpt{d_{{\rm pt}}}
\def\sjz{\sum_{j\in\integ}}
\def\prax{Proposition~\ref{appx}}

\def\uii{^{(i)}}
\def\twi{{\widetilde W}\uii}
\def\uil{u_{il}}
\def\vil{v_{il}}
\def\sln{\sum_{l=1}^n}
\def\BP{{\rm BP}}

\title{On Stein's method and perturbations}
\author
{A.\ D.\ Barbour\thanks{Angewandte Mathematik,
Winterthurerstrasse~190,
CH--8057 Z\"URICH, Switzerland: {\tt a.d.barbour@math.unizh.ch}
},\, Vydas  \Ceka\thanks{Department of Mathematics and 
Informatics, Vilnius University,
Naugarduko 24, Vilnius 03225, Lithuania: {\tt
vydas.cekanavicius@mif.vu.lt}}~ and Aihua
Xia\thanks{Department of Mathematics and Statistics, University of
Melbourne, Vic 3010, Australia: {\tt xia@ms.unimelb.edu.au}}
\\
Universit\"at Z\"urich, Vilnius University\\ and University of Melbourne}
\date{}
\maketitle

\vglue-1cm
\begin{abstract}
Stein's~(1972) method is a very general tool for assessing 
the quality of approximation of the distribution of a random element 
by another, often simpler, distribution. In applications of Stein's 
method, one needs to establish a Stein identity for the approximating 
distribution, solve the Stein equation and estimate the behaviour of 
the solutions in terms of the metrics under study. For some Stein 
equations, solutions with good properties are known; for others, this 
is not the case.  
\BX~(1999) introduced a perturbation method 
for Poisson approximation, in which Stein 
identities for a large class of compound Poisson and translated Poisson 
distributions are viewed
as perturbations of a Poisson distribution. In this paper, 
it is shown that the method can be extended to very general 
settings, including perturbations of normal, Poisson, compound Poisson, 
binomial and Poisson process approximations in terms of various metrics 
such as the Kolmogorov, Wasserstein and total variation metrics. Examples 
are provided to illustrate how the general perturbation method can 
be applied.
\end{abstract}

\vskip10pt

\nin {{\it Keywords:} perturbation method, normal distribution, 
jump diffusion process, Poisson distribution, compound Poisson 
distribution, Poisson process, point process, total variation norm, 
Kolmogorov distance, Wasserstein distance, local distance.}

\section{Introduction}\label{intro}
 \setcounter{equation}{0}
Many applications of Stein's~(1972) method, when approximating the
distribution~$\law(W)$ of a random element~$W$ of a metric space~$\xx$
by a probability distribution~$\pi$, are accomplished broadly as
follows.  The aim is to estimate 
$\ex h(W)-\pi(h)$ for each member~$h$ of a family of test functions $\hh$, 
where $\pi(h):=\int hd\pi$. To do this, one finds a normed space~$\gg$ and 
an appropriate Stein operator~$\aaa$ on $\gg$ characterizing~$\pi$; 
$\aaa\colon \gg \to \ff \subset \re^{\xx}$, for some $\ff\supset\hh$, 
must be such that $\pi(\aaa g) = 0$ for all~$g$ in~$\gg$, and 
that~$\pi$ is the unique probability distribution for which this is
the case. `Appropriate' in this context means that an inequality of
the form
\eq\label{intro-1}
  |\ex\{(\aaa g)(W)\}| \Le \e\|g\|_{\gg},\qquad g \in \gg,
\en
can be established, for some (small)~$\e$.  
%\adb{Stein's method for~$\pi$
%based on~\Ref{intro-1} can then be said
%to `work' if~\Ref{intro-1} is already enough to imply that $\law(W)$
%and~$\pi$ are close, usually in the sense that the distance between them,
%measured with respect to some chosen distance, is small with~$\e$.}
%
%\adb{For instance, suppose that $\hh\subset\ff$
%is a family of test functions and that, for each $h\in\hh$, there is
Finally, for each $h\in\hh$, find 
a function $g_h \in\gg$ satisfying the Stein equation
\eq\label{intro-2a}
  \aaa g_h \Eq h - \pi(h).  
\en
Then it follows from~\Ref{intro-1} that
\eq\label{intro-2}
  |\ex h(W) - \pi(h)| \Le \e\|g_h\|_\gg.
\en
Hence, if it can be shown that 
\eq\label{intro-3}
   \|g_h\|_\gg \Le C\|h\|_\ff,
\en
for some norm $\|\cdot\|_\ff$ on~$\ff$, we can conclude that
\eq\label{intro-3a}
  d_\hh(\law(W),\pi) \Le C\e \sup_{h\in\hh} \|h\|_\ff,
\en
where, for any two distributions $P$ and $Q$ on $\xx$,
\eq\label{intro-4}
   d_\hh(P,Q) \ :=\ \sup_{h\in\hh} |P(h)-Q(h)|.
\en
Thus, if \Ref{intro-2a} and~\Ref{intro-3} are satisfied, it
is enough for the $d_\hh$-approximation of $\law(W)$ by~$\pi$ 
to establish the inequality~\Ref{intro-1};  in this sense, Stein's
method for~$\pi$ can be said to work for the distance~$d_\hh$.
Distances of this form include the total variation distance $d_{TV}$, with
$\hh$ the set of functions bounded by~$1$, and the Wasserstein
distance $d_W$, with~$\hh$ the Lip\-schitz functions with slope 
bounded by~$1$.

\adb{
Probabilistic inequalities of the form~\Ref{intro-1} can be derived
by a variety of techniques, including Stein's exchangeable pair
approach, the generator method and Taylor expansion.  However, the
analytic inequality~\Ref{intro-3} can prove to be a stumbling block,
especially if a reasonably small value of~$C$ is desired, unless~$\pi$
happens to be a particularly convenient distribution. For $\xx = \re$,
the normal and Poisson distributions lead to simple versions 
of~\Ref{intro-3}. However, when introducing Stein's method for compound 
Poisson distributions, \BCL~(1992) were only able to prove analogous
inequalities with satisfactory values of~$C$ for distributions for which
the generator method was applicable, and this represents a strong 
restriction on the compound Poisson family.  The class of amenable
compound Poisson distributions was subsequently extended in \BX~(1999),
where a perturbation technique was introduced, which enabled
the good properties 
of the solutions of the Poisson operator to be carried over to those
of the Stein equations for neighbouring compound Poisson distributions.
Their approach was taken further in
Barbour \&  \Ceka~(2002) and in \Ceka~(2004). Here, we show that
the perturbation idea can be applied not just in the Poisson setting,
but in great generality. One consequence is that
the range of compound Poisson distributions whose solutions
have good properties can be further extended, but the scope of possible
applications is much wider. In particular, there is no need to restrict
attention to random variables on the real line; distributions and random 
elements on quite general spaces can be considered.}

\adb{
The perturbation method is discussed in the general terms in Section~\ref{formal}.
Theorem~\ref{perturbation} shows how to find the solution~$g_h$ in~\Ref{intro-2a}
for $\aaa = \aaa_1$, when~$\aaa_1$ is close enough to a `nice' Stein
operator~$\aaa_0$, and the probability measure~$\pi_0$ associated with~$\aaa_0$
has $\supp(\pi_0) = \xx$;  the theorem also gives the inequality
corresponding to~\Ref{intro-3}.  Theorem~\ref{very-useful} gives conditions
under which Stein's method works, but which do not assume the support
condition, and Theorem~\ref{K-distance} allows a further slight relaxation,
which is particularly relevant to approximation of random variables
using the Kolmogorov distance.}
In Section~\ref{examples}, a number of specific examples are given, some of
which are illustrated from the point of view of application in 
Section~\ref{illustration}.

\adb{
As indicated above, there are various ways in which an 
inequality~\Ref{intro-1} relevant in any particular setting may be derived.
This means that the choice of operator~$\aaa_1$, and of the corresponding
approximating probability measure~$\pi_1$, is frequently dictated by the
problem under consideration in a more or less natural way.  The choice
of~$\aaa_0$ is more a matter of chance.  If~$\aaa_1$ is not 
itself one of the operators for which the solutions to~\Ref{intro-2a}
are known to satisfy an inequality of the form~\Ref{intro-3}, then one
looks for an~$\aaa_0$ which is, and which is not too far away from~$\aaa_1$.
Such an operator need not exist.  In order for our perturbation approach
to be successful, it is necessary for the contraction inequality~\Ref{key}
to be satisfied, and this limits the set of operators which can be considered
as perturbations of any given~$\aaa_0$, for the purposes of our theorems.}

\section{Formal approach}\label{formal}
\setcounter{equation}{0}
Let $\xx$ be a Polish space, and~$\gg$ a linear subspace of the functions 
$g\colon \xx\to\re$ equipped with a norm $\|\cdot\ngg$.
Suppose that~$\pi_0$ is a probability measure on~$\xx$ with $\supp(\pi_0) = \xx_0
\subset \xx$.   Define
\eqs
  \ff &:=& \{f\colon \xx\to\re,\ \pi_0(|f|) < \infty\};\\
  \ff_0 &:=& \{f\in\ff\colon  f(x) = 0\ 
      \mbox{for all}\ x \notin \xx_0\};\\
  \ff' &:=& \{f\in\ff\colon \pi_0(f)=0\}; \qquad \ff'_0 \ :=\ \ff_0 \cap \ff',
\ens
and let $P_0$ be the projection from $\ff$ onto $\ff'_0$ given by
\[
  P_0f \ :=\ f\ixo - \pi_0(f)\ixo,
\]
where, here and subsequently, $\bone_A$ denotes the indicator function 
of the set~$A$, and multiplication of
functions is to be understood pointwise.  Now let $\|\cdot\|$ be a norm on~$\ff$, 
set
\[
  \bff \ :=\ \{f\in\ff\colon \|f\| < \infty\},
\]
and define $\bff_0 := \bff\cap\ff_0$, $\bff' := \bff\cap\ff'$,
$\bff'_0 := \bff\cap\ff'_0$; we shall require that $\|\cdot\|$ is such that 
\eq\label{norm-cond1}
  P_0\colon \bff \to \bff'_0.
\en
\ignore{
and that, for any $f\in\bff$,
\eq\label{norm-cond3}
  |\pi_0(f)| \Le \k_0 \|f\|.
\en
These conditions are clearly always satisfied with $\k_0=1$
for the supremum norm~$\|\cdot\nsup$; for an $L_p$-norm with respect to a measure~$\mu$
on~$\xx$, $p \ge 1$, they are satisfied if~$\pi_0$ has bounded density
with respect to~$\mu$.
}  
We also assume that~$\bff$ is a determining class of functions for
probability measures on~$\xx$ (Billingsley~1968, p.~15). 

We now suppose that there is a `nice' Stein operator $\aaa_0$ characterizing~$\pi_0$.
By this, we mean that
\eq\label{Stein-operator0-1}
  \aaa_0\colon \gg \to \ff_0', 
  %\quad\mbox{and}\quad (\aaa_0 g)(x) = 0\ \mbox{for all}
  %\ x\notin \xx_0,
\en
and also that it is possible to define a right inverse
\[
  \aaomi\colon \bff'_0\ \to\ \gg_0\ :=\ \{g\in\gg\colon g(x) = 0\ 
      \mbox{for all}\ x \notin \xx_0\},
\]      
satisfying
\eqa
   \aaa_0(\aaomi f) \Eq f \quad \mbox{for all} \quad f \in \bff'_0;\label{Stein-operator0-2}\\
   \|\aaomi P_0 f\ngg \Le A\|f\|, \qquad f \in \bff, \label{Stein-operator0-3}
\ena
for some $A < \infty$. Note that~\Ref{Stein-operator0-1} means that
\eq\label{characterization}
  \pi_0(\aaa_0 g) = 0\quad \mbox{for all}\quad g\in\gg. 
\en
On the other hand, in view 
of~\Ref{Stein-operator0-2}, if~$\pi$ is any probability measure on~$\xx_0$ such
that $\pi(\aaa_0 g) = 0$ for all $g\in\gg$, then $\pi(f) = 0$ for all $f \in \bff'_0$,
meaning that $\pi(f) = \pi_0(f)$ for all $f\in\bff_0$,
and hence for all $f\in\bff$. Since~$\bff$ is a
determining class, $\pi = \pi_0$, and~$\aaa_0$ characterizes~$\pi_0$ 
through~\Ref{characterization}.

\adb{
In the setting of the introduction, for $h \in \hh \subset \ff_0$ a family of
test functions, we have $h(x) - \pi_0(h) = (P_0h)(x)$ for $x \in \xx_0$, so
that we can take $g_h = \aaa_0^{-1}P_0h$ and obtain~\Ref{intro-2a}, in view
of~\Ref{Stein-operator0-2}.
Inequality~\Ref{Stein-operator0-3} is just~\Ref{intro-3}
for~$\aaa_0$, with $f$ in place of~$h$.  Hence, because of~\Ref{intro-3a}, 
Stein's method for~$\pi_0$ based
on~\Ref{intro-1} (with $\aaa_0$ in place of~$\aaa$) works for distances based
on families~$\hh$ of test functions whose norms are uniformly bounded.
  Our interest here is in extending this to
probability measures~$\pi_1$ characterized by generators~$\aaa_1$ which are
close to~$\aaa_0$.}

So let~$\pi_1$ be a finite signed measure on~$\xx$ with $\pi_1(\xx) = 1$, 
and such that $|\pi_1|(|f|) < \infty$ for all $f\in\bff$.
\ignore{  
and that, for some fixed~$\k < \infty$,
\eq\label{norm-cond2}
  |\pi_1(f\ixoc)| \Le \k_1|\pi_1|(\xx_0^c)\,\|f\|,\qquad f\in\bff.
\en
Once again, the latter condition is always satisfied with $\k_1=1$
for the supremum norm, and also for an $(L_p,\mu)$-norm, $p \ge 1$, 
if~$d\pi_1/d\mu$ is bounded on~$\xx_0^c$.
} 
Let~$\aaa_1$ be a
Stein operator for~$\pi_1$, meaning that $\aaa_1\colon \gg \to \ff'_1$, where
\[
  \ff'_1 \ :=\ \{f\colon \xx\to\re;\ |\pi_1|(|f|) < \infty,\ \pi_1(f) = 0\},
\]
so that
\eq\label{Stein-operator1-1}
  \pi_1(\aaa_1 g) \Eq 0 \quad \mbox{for all} \quad g \in \gg;
\en
set $\uu = \aaa_1 -\aaa_0$, and assume also that 
\eq\label{Stein-operator1-2}
%  \aaa_1(\aaomi P_0)\colon \bff \to \bff \quad \mbox{and}\quad 
  \uu\aaomi P_0\colon \bff\to\bff.
\en
The key assumption which ensures that~$\aaa_1$ can fruitfully be thought of as
a perturbation of~$\aaa_0$ is that
\eq\label{key}
  \|\uu\aaomi P_0\| \ =:\ \g \ <\ 1.
\en

\nin{\bf Remark.} Having to satisfy the condition~\Ref{key} significantly limits the choice of
distributions~$\pi_1$ whose Stein equations can be treated as perturbations
of that for~$\pi_0$.  This is clearly illustrated in the examples of the next
section.

\begin{theorem}\label{perturbation}
With the above definitions, suppose that assumptions 
\Ref{norm-cond1}--\Ref{Stein-operator0-3} and \Ref{Stein-operator1-1}--\Ref{key}
are satisfied. Then the operator
\eq\label{B-def}
  \bb \ :=\ \aaomi P_0 \sjo (-1)^j(\uu\aaomi P_0)^j\colon\  \bff\ \to\ \gg_0
\en
is well defined, and 
\eq\label{B-norms}
  \|\bb\| \Le A/(1-\g);\qquad \|\uu\bb\| \Le \g/(1-\g).
\en
Furthermore, for $f\in\bff$ and for all  $x \in \xx_0$,
\eq\label{A1-solution}
  (\aaa_1\bb f)(x) - (P_1f)(x) \Eq c(f) \Eq \pi_1(f) - \pi_0(f) + \pi_0(\uu\bb f),
\en
where $P_1 f = f - \pi_1(f)\bone$; here, $\bone=\bone_\xx$.
\ignore{
, and 
\eq\label{c-bnd}
  |c(f)| \Le  \frac{2|\pi_1|(\xx_0^c)}{1-\g}\,\|f\nsup.
\en
}
In particular, if $\xx_0=\xx$, we have $c(f)=0$, so that~$\bb$ is a right
inverse of~$\aaa_1$ on~$\ff'_1\cap\bff$.
\end{theorem}
         
\begin{proof}
The first part is immediate from \Ref{Stein-operator0-3} and \Ref{key}, 
from the properties of $\aaomi$
and from~\Ref{Stein-operator1-2}.  It is then also immediate that
\[
  (\aaa_0 + P_0\uu)\bb f \Eq P_0 f,\quad f\in\bff.
\]
Hence, for $f\in\bff$, we have
\eqa
  \aaa_1\bb f &=& (\aaa_0 + P_0\uu + (I-P_0)\uu)\bb f \non\\
  &=& P_0f + (\uu\bb f)\ixoc + \pi_0(\uu\bb f)\ixo, \label{a1bf-1}
\ena
so that, for $x\in\xx_0$,
\eq\label{a1bf-2}
  (\aaa_1\bb f)(x) - (P_1 f)(x) \Eq  \pi_1(f) - \pi_0(f) + \pi_0(\uu\bb f)
     \ =:\ c(f).
\en 
      
For the constant~$c(f)$, note that, from~\Ref{Stein-operator1-1} with
$\bb f$ for~$g$ and from~\Ref{a1bf-1}, we have
\eqs
  0 &=& \pi_1(f\ixo) - \pi_1(\xx_0)\pi_0(f)
      + \pi_1((\uu\bb f)\ixoc) + \pi_0(\uu\bb f)\pi_1(\xx_0)\\
  &=& \pi_1(f) - \pi_1(f\ixoc) - \pi_0(f) + \pi_1(\xx_0^c)\pi_0(f)
     + \pi_1((\uu\bb f)\ixoc)\\
  &&\hspace{2cm}\mbox{}  + \pi_0(\uu\bb f)(1 - \pi_1(\xx_0^c)).
\ens
This implies, 
% from \Ref{norm-cond3} and~\Ref{norm-cond2} and 
from the first part of the theorem, that
\eqs
  c(f) &=&  \pi_1(f) - \pi_0(f) + \pi_0(\uu\bb f) \Eq 0
\ens
if $\xx_0^c = \emptyset$, and
\eq\label{c(f)-formula}
   c(f)  \Eq  \pi_1(f\ixoc) -  \pi_1(\xx_0^c)\pi_0(f) - \pi_1((\uu\bb f)\ixoc)
       +  \pi_0(\uu\bb f)\pi_1(\xx_0^c) 
\en
otherwise. \ep
%  &\le&   |\pi_1|(\xx_0^c)\,\|f\nsup(1 + 1 +  \{\g/(1-\g)\}  + 
 %      \{\g/(1-\g)\}),
% \ens
% as required.\ep
\end{proof}

\bigskip
\nin{\bf Remark.}
\adb{If $\xx_0 = \xx$, then it follows from Theorem~\ref{perturbation} that
\[
   \aaa_1\bb h \Eq P_1h \Eq h - \pi_1(h)
\]
for test functions $h \in \hh \subset \bff$.  Hence, for such~$h$,
the function $g_h := \bb h$ satisfies~\Ref{intro-2a}, where $\aaa$
is replaced by~$\aaa_1$ and $\pi$ by~$\pi_1$. It then follows
from~\Ref{B-norms} that $\|g_h\ngg \le A(1-\g)^{-1}\|h\|$, so 
that~\Ref{intro-3} is satisfied with $C = A/(1-\g)$, and hence
Stein's method for~$\pi_1$ based on~\Ref{intro-1} (with $\aaa_1$ for~$\aaa$)
works for distances~$d_\hh$ derived from bounded families of test functions.}

\adb{If $\xx_0 \neq \xx$, the inequalities~\Ref{B-norms} are still satisfied,
so that~\Ref{intro-3} is still true with $C = A/(1-\g)$ if $g_h = \bb h$.
However, this choice of~$g_h$ now gives only an approximate solution 
to~\Ref{intro-2a}:
\eq\label{approx-2a}
  (\aaa_1 g_h)(x) \Eq h(x) - \pi_1(h) + c(h), \qquad x \in \xx_0.
\en
This is still enough to show that Stein's method works for~$\pi_1$ based 
on~\Ref{intro-1} (with $\aaa_1$ for~$\aaa$), as is demonstrated
in Theorem~\ref{very-useful} below.  To make the connection, we first
need two more lemmas.}		 

\bigskip
\adb{The first concerns the size of $|c(f)|$.  This can be controlled
in a number of ways,  two
of which are} given in the following lemma.  For any finite signed 
measure~$\pi$ and any $A \subset \xx$, we define
\eq\label{kappa-def}
  \k(\pi,A) \ :=\ \sup_{\{f\in\bff\colon \|f\| \le 1\}} |\pi| (\hat f\bone_A),
\en
\adb{where $\hat f(x) := |f(x) - \pi_0(f)|$.}

\begin{lemma}\label{c(f)-bnds}
For $f\in\bff$, we have
\eqs
  \mbox{(i)}&& |c(f)| \Le  \frac{2|\pi_1|(\xx_0^c)}{1-\g}\,\|f\nsup;\\
 \mbox{(ii)}&& |c(f)| \Le \frac{\k(\pi_1,\xx_0^c) }{1-\g}\,\|f\|.
\ens
% where $\k_0 := \k(\pi_0,\xx_0)$.
% \sup_{\{f\in\bff\colon \|f\| \le 1\}}\pi_0(|f|)$.
\end{lemma}

\begin{proof}
The proof is immediate from~\Ref{c(f)-formula} and~\Ref{kappa-def}.
\ep  \end{proof}
	         
\bigskip
\adb{The second lemma translates~\Ref{approx-2a} into an inequality bounding
the difference
$|\pi(f) - \pi_1(f)|$ in terms of $|\pi(\aaa_1\bb f)|$, for a general
probability measure~$\pi$ on~$\xx$.}

\begin{lemma}\label{useful}
Under the conditions of Theorem~\ref{perturbation}, if~$\pi$ is any probability
measure on~$\xx$, then, for any $f\in\bff$, we have
\eqs
    |\pi(f) - \pi_1(f)| \Le |\pi(\aaa_1\bb f)| 
%\\ &&\quad\qquad\mbox{}
 +  \begin{cases}
     2(1-\g)^{-1}\{|\pi_1|(\xx_0^c) + \pi(\xx_0^c)\}\,\|f\nsup;\\
    (1-\g)^{-1} \{\k(\pi_1,\xx_0^c) + \k(\pi,\xx_0^c)\}\,\|f\|.
  \end{cases}
\ens
\end{lemma}

\begin{proof}
It follows from~\Ref{a1bf-1} that
\eqs
  \pi(\aaa_1\bb f) &=& \pi(P_0 f) + \pi((\uu\bb f)\ixoc) + \pi_0(\uu\bb f)\pi(\xx_0)\\
  &=& \pi(f) - \pi(f\ixoc) - \pi_0(f)(1 - \pi(\xx_0^c)) \\
  &&\hspace{2cm}\mbox{} + \pi((\uu\bb f)\ixoc) + \pi_0(\uu\bb f)(1-\pi(\xx_0^c))\\
  &=& \{\pi(f)-\pi_1(f)\} + c(f) - \pi(f\ixoc) \\
  &&\hspace{2cm}\mbox{} + (\pi_0(f) - \pi_0(\uu\bb f))\pi(\xx_0^c)
    + \pi((\uu\bb f)\ixoc).
\ens
Hence, and using~\Ref{kappa-def}, the lemma follows.\ep
\end{proof} 

\bigskip
\adb{This lemma gives the information that we need,} when deriving 
distributional approximations
in terms of the measure~$\pi_1$.  Let~$\hh \subset \bff$ be any collection
of test functions which 
%is closed under scalar multiplication and
forms a determining class for probability measures on~$\xx$.
Then define the metric~$\dhh$ on finite signed
measures $\r,\s$ on~$\xx$, by
\eq\label{H-dist}
  \dhh(\r,\s) \ :=\ \sup_{h\in\hh} |\r(h) - \s(h)|.
\en  
In the special case where $\hh := \{f\in\bff\colon \|f\| \le 1\}$, 
we write $\dff$ for~$\dhh$.  \adb{The following theorem shows that
Stein's method for~$\pi_1$ based on~\Ref{Stein-approx} works
for the distance~$\dff$, even when $\xx_0 \neq \xx$.}

\begin{theorem}\label{very-useful}
\adb{Suppose that the conditions of Theorem~\ref{perturbation} are
satisfied, and write $g_f^0 := \aaomi P_0 f$ for all $f\in\bff$.
Then, if 
\eq\label{Stein-approx}
 | \pi(\aaa_1 g_f^0)| \Le \e\|g_f^0\ngg \quad \mbox{for all}\quad f\in \bff,
\en
it follows that
\[
  \dff(\pi,\pi_1) 
    \Le (1-\g)^{-1}\{A\e + \e'(\pi,\pi_1)\}, 
\]
where}
\eqs
  \e'(\pi,\pi_1) \ :=\  \min\{2(|\pi_1|(\xx_0^c) + \pi(\xx_0^c))F, 
      \k(\pi_1,\xx_0^c) + \k(\pi,\xx_0^c) \},
\ens
and
\eq\label{F-def}
     F := \sup_{\{f\in\bff\colon \|f\| \le 1\}}\|f\nsup.
\en 
\end{theorem}

\begin{proof} In fact, let $\tilde f=\sum_{j\ge 0}(-1)^j(\uu\aaomi P_0)^jf$.
Then~\Ref{Stein-approx} \adb{together with~\Ref{key} and~\Ref{Stein-operator0-3} 
imply that
\[
  |\pi(\aaa_1\bb f)|=|\pi(\aaa_1 g_{\tilde f}^0)| \Le  \e\|g_{\tilde f}^0\ngg
	   \Le A\e\|\tilde f\| \Le \frac{A\e}{1-\g}\|f\|.
\]
Thus the conclusion follows immediately from Lemma~\ref{useful}
and} from the definition of~$\dff$.
\ep
\end{proof}

\bigskip
\adb{\nin Note that~\Ref{Stein-approx} is a weakening of what would
normally be required for~\Ref{intro-1}, inasmuch as the inequality
is only needed for the functions~$g_f^0$, which, being the
solutions to the Stein equation for the `nice' operator~$\aaa_0$,
may well be known in advance to have good properties.}

Theorem~\ref{very-useful} is applied most simply when~$\pi$ is the 
distribution of some random
element~$W$, for which it can be shown that
\eq\label{error}
  |\ex(\aaa_1 g)(W)| \Le \sum_{j=1}^l \e_j c_j(g), \quad g\in\gg.
\en       
Here, the quantities~$\e_j$ are to be computed using~$W$ alone, and 
the function~$g$ enters
only through the constants~$c_j(g)$.  If the norm~$\|\cdot\|$ on~$\ff$ can be
chosen in such a way that the $c_j(g_f^0)$ can be bounded by a multiple of~$\|f\|$
for any $f\in\bff$, then Theorem~\ref{very-useful} can be invoked.

The choice of norms on~$\ff$ for which this procedure can be carried through depends
very much on the structure of the random variable~$W$: see Section~\ref{illustration}.
Broadly speaking, for the more stringent norms, the contraction
condition~\Ref{key} is harder to satisfy;
on the other hand, there are then fewer functions having finite norm,
and so the inequality~\Ref{Stein-approx} is easier
to establish.  Take, for example,
% The need for choices of~$\hh$ smaller than the whole of~$\bff$
% is illustrated already in the case of 
standard normal approximation, with~$\gg$ the space of bounded real functions
with bounded \adb{first and second derivatives, endowed with the norm
\eq\label{normal-G}
  \|g\ngg \ :=\ \|g\nsup + \|g'\nsup + \|g''\nsup \, ,
\en
and} with~$\aaa_0$ the Stein operator given by
\eq\label{Stein-standard}
  (\aaa_0 g)(x) \Eq g'(x) - x g(x),    \quad g\in \gg.
\en
Here, it is possible, in many central limit settings, to derive an
inequality of the form~\Ref{intro-1}: 
\[
   |\ex(\aaa_0 g)(W)| \le  \e\|g\ngg
\]
for some~$\e$, as, for example, in Chen \& Shao~(2005, p.~5). 
% given in~\Ref{error}, with $\aaa_0$ for~$\aaa_1$,
% $l=1$ and $c_1(g) \le k\|g''\nsup$ for some~$k$ **tutorial reference**. 
Now, for $g_f^0=\aaomi P_0 f$ with 
$\|f'\nsup < \infty$, we \adb{have $\|(g_f^0)''\nsup \le 4\|f'\nsup$
by \prax\,(c)(i) and~(iii) with $y=g_f^0$, 
so that inequality~\Ref{intro-3} is satisfied with
\eq\label{diff-norm}
   \|f\|\ui \ :=\  \|f\nsup + \|f'\nsup
\en
as norm on~$\ff$.
This, in turn, leads to corresponding approximations with respect to the
distance~$\dff =: d\ui$, from~\Ref{intro-3a}.}

In the usual central limit context, there is typically no  hope of
taking the argument further, and choosing $\hh = \bff$ for the
supremum norm~$\|\cdot\nsup$ in place of~$\|\cdot\|\ui$ on~$\ff$.  
This is not because the perturbation
argument would fail, but because there can usually be no inequality
of the form $|\ex(P_0 f)(W)| \le \e\|f\nsup$ for all $f\in\bff$, 
unless~$\e$ is rather large; this is because the supremum of
the left hand side is then just the {\it total variation\/} distance between
$\law(W)$ and the standard normal distribution, and this is not
necessarily small under the usual conditions for the central limit
theorem. More is, however, possible with some extra restrictions: 
see Cacoullos \etal~(1994) and Example~\ref{illustration}.1. 

The distance~$d\ui$ is not the one most commonly used for
measuring the accuracy of approximation in the central limit theorem.
Here, it is usual to work with the Kolmogorov distance~$d_K$, which is
of the form defined in~\Ref{H-dist}, with the set of
test functions
\[
  \hh^K\ :=\ \{\bone_{(-\infty,a]}\colon a \in \re\}.
\]
For these test functions, it can in many central limit
applications be established, albeit with rather more effort,
that $|\ex(\aaa_0 g_h)(W)|$ is bounded, \adb{uniformly for $h \in \hh^K$,} 
by a quantity of the form 
$k\e$ for some $k<\infty$ and~$\e$ reflecting  the
closeness of~$\law(W)$ and the standard normal distribution.
This in turn, with~\Ref{intro-2a}, implies  error estimates for standard normal
approximation, measured with respect to Kolmogorov distance.

Now the set~$\hh^K$ forms a subset of~$\bff$, when the supremum norm is
taken on~$\ff$, and  the
perturbation arguments leading to Lemma~\ref{useful} can still
be applied successfully, for Stein operators~$\aaa_1$ suitably close
to~$\aaa_0$.  However, in order to deduce distance estimates
as in Theorem~\ref{very-useful}, it is necessary to be able to
bound $|\ex(\aaa_0 g_f^0)(W)|$ \adb{not only for $f \in \hh^K$, but also 
for any~$f$ of the form $f := (\uu\aaomi P_0)^j h$,
where $h\in\hh^K$ and $j\ge1$, since these functions are used to make 
up the function~$\tilde f$ introduced in the proof of Theorem~\ref{very-useful}. Now these 
functions~$f$ are not typically
in the set~$\hh^K$. However, it can at least be shown that both
$g_h^0$ and~$(g_h^0)'$ are uniformly bounded  for $h\in\hh^K$. For some
operators~$\aaa_1$, this is enough to be able to conclude that
\[
   \sup_{h\in\hh^K} \|\uu g_h\|\ui \ <\ \infty.
\]
It is then possible to apply} the following result, in which the Stein operator~$\aaa_0$ 
is now quite general.

\begin{theorem}\label{K-distance}
  Suppose that the conditions of Theorem~\ref{perturbation} are satisfied,
and that~$\hh$ is any family of test functions with 
$H := \sup_{h\in\hh}\|h\nsup < \infty$, and such that $g_h^0 := \aaomi P_0 h$ 
is well defined for $h\in\hh$, satisfying $\aaa_0 g_h^0 = P_0 h$ and $\uu g_h^0 \in \bff$.  
Assume further that
\eq\label{U-cond}
  \g_\hh \ :=\ H^{-1}\sup_{h\in\hh}  \|\uu g_h^0\| \ <\ \infty.
\en
% where $g_f^0 := \aaomi P_0 f$ for $f\in\bff$.  
Then, if~$\pi$ is such that
\eq\label{cond-1}
  \sup_{h\in\hh} |\pi(\aaa_1 g_h^0)| \Le H\e_1
\en
and
\eq\label{cond-2}
  |\pi(\aaa_1 g_f^0)| \Le \e_2\|g_f^0\ngg,\quad f \in \bff,
\en
it follows that
\[
  \dhh(\pi,\pi_1) \Le H\Blb \e_1 + \frac{\g_\hh A\e_2}{1-\g} + 
    \frac{\e(\pi,\pi_1)}{1-\g}\Brb,
\]
where $\e(\pi,\pi_1):=\k(\pi_1,\xx_0^c) + \k(\pi,\xx_0^c)$.
\end{theorem}

\begin{proof}
Once again, \adb{much as in the proof of Theorem~\ref{very-useful}, we
note that
\eq\label{K-estimate}
  |\pi(\aaa_1\bb h)| \Le \sjo|\pi(\aaa_1\aaa_0^{-1}P_0(\uu\aaa_0^{-1}P_0)^j\, h)|
	  \Eq |\pi(\aaa_1 g_h^0)| + \sji|\pi(\aaa_1 g_{f_j}^0)|,
\en
where $f_j :=(\uu\aaa_0^{-1}P_0)^j\,h$, $j\ge1$.  Now, for $h\in\hh^K$,
\[
   \|f_1\| \Eq \|\uu g_h^0\| \Le H\g_H,
\]
by \Ref{U-cond}, and then, by~\Ref{key},
\[
  \|f_j\| \Le \g^{j-1} H\g_H,\qquad j\ge2.
\]
Hence, from~\Ref{K-estimate}, \Ref{Stein-operator0-3}, \Ref{cond-1}
and~\Ref{cond-2}, it follows that
\[
   \sup_{h\in\hh^K} |\pi(\aaa_1\bb h)| \Le H\e_1 + \sji \e_2 A\g^{j-1} H\g_H,
\]
and the theorem now follows from Lemma~\ref{useful}.}	 		  		
\ep
\end{proof}

\bigskip
In particular, if~$\aaa_0$ is the Stein operator for normal approximation
given in~\Ref{Stein-standard}, and taking the norm~$\|\cdot\|\ui$,
Theorem~\ref{K-distance} can be applied with $\hh = \hh^K$; in
circumstances in which the conditions \Ref{U-cond}--\Ref{cond-2}
are satisfied, this leads to estimates of the error in approximating
the distribution~$\pi$ of a random variable~$W$ by~$\pi_1$,
measured with respect to Kolmogorov distance.  In particular,
the estimates \Ref{cond-1} and~\Ref{cond-2} relating to the distribution
of~$W$ are of a kind which can often be verified in practice; see
Section~\ref{illustration}.

\ignore{ satisfied 
$\hh'$ is not usually rich enough
to fulfil the conditions of Theorem~\ref{very-useful}, in that one
cannot expect to have $\uu\aaomi P_0 \hh \subset \hh$.  To get
round this difficulty, one can, for instance, take for~$\hh$ the
larger collection of functions
\eq\label{HK-def}
  \hh^K \ :=\ \{h_1+ch_2\colon h_1\in\bff\ui,\, h_2\in\hh',\,c\in\re\}
  \supset \hh'.
\en
For $h\in\hh^K$ and for~$\aaa_0$ given in~\ref{Stein-standard}, it
follows that $\aaomi P_0 h \in \bff\ui$, and this is sometimes
enough to conclude that $\uu\aaomi P_0 \hh \subset \hh$ ---
see, for example, Examples \ref{examples}.1 and~\ref{examples}.2
below.  What is more, arguments used to show that $\ex\{(\aaa_1 g_h)(W)\}$ 
is uniformly small for $h\in\hh'$ can be expected to apply also to
functions $h\in\hh^K$, since the $\bff\ui$-component of such an~$h$
is treated easily as discussed around~\Ref{diff-norm}.  In such
circumstances,  Theorem~\ref{very-useful} can be
invoked to deduce error estimates with respect to Kolmogorov
distance when approximating by perturbations~$\pi_1$ of the normal
distribution, provided  that the conditions
of Theorem~\ref{perturbation} are satisfied.
}

\ignore{
In particular, applying~$\pi_0$ to~\Ref{A1-solution}
for any $f\in\bff$ gives
\[
   \pi_0(\aaa_1\bb f) - \pi_0(f) + \pi_1(f) \Eq c(f),
\]
and, from~\Ref{a1bf-1}, 
\[
  \pi_0(\aaa_1\bb f) \Eq \pi_0(\uu\bb f);
\]
hence, from \Ref{B-norms} and Lemma~\ref{c(f)-bnds}, $\pi_1$ is constrained to be fairly 
close to~$\pi_0$, in the sense that
\eq\label{pi0-pi1}
  d_{\bff}(\pi_0,\pi_1) \Le  \frac{\k_0\g + 2F|\pi_1|(\xx_0^c)}{1-\g},
\en
where~$F$ is as in~\Ref{F-def}.
}
 
\section{Examples}\label{examples}
\setcounter{equation}{0} 
In the first two examples, the sets $\xx_0$ and~$\xx$ are the same, so that the
elements in the bounds involving probabilities of the set~$\xx_0^c$ make no
contribution.  The first of these is purely for illustration, since properties
of the Stein equation for the perturbed distribution could be
obtained directly.

\medskip     
\nin{\bf Example~\ref{examples}.1}.\ 
In this example, we consider approximation by the probability 
distribution~$\pi_1 := t_{m,\psi}$ on~$\re$, with density
\[
    p_{m,\psi}(x) \Eq  k_{m,\psi}(1+x^2/m)^{-(m+1)\psi/2}\,e^{-(1-\psi)x^2/2},
      \quad x\in\xx := \re,
\]
where~$k_{m,\psi}$ is an appropriate normalizing constant.  
\adb{ This family of densities interpolates between the standard normal 
($\ps=0$) and Student's~$t_m$ distribution ($\ps=1$)
distribution, as~$\ps$ moves from $0$ to~$1$; $m$ is classically a
positive integer.}  We take for~$\gg$ the space of bounded real functions
with bounded derivatives, endowed with the norm
\[
  \|g\ngg \ :=\ \|g\nsup + \|g'\nsup\,.
\]
An appropriate Stein operator~$\aaa_1$ for $t_{m,\psi}$ is given by
\eq\label{t-Stein}
  (\aaa_1 g)(x) \Eq g'(x) - x\Blb (1-\psi) + \frac{\psi(m+1)}{m+x^2} \Brb g(x),
    \quad g\in \gg;
\en
\adb{this follows because $p_{m,\ps}(x)$ is an integrating factor for the
right hand side of~\Ref{t-Stein}, and hence, for any $g\in\gg$,
\[
  \int_{-\infty}^\infty (\aaa_1 g)(x) p_{m,\ps}(x)\,dx 
	  \Eq \left[g(x)p_{m,\ps}(x) \right]_{-\infty}^\infty \Eq 0,
\]
so that~\Ref{Stein-operator1-1} is satisfied.		
Now, at least for small enough~$\psi$, $\aaa_1$ could be} thought of as a 
perturbation of the standard normal distribution, with Stein operator
\[ 
   (\aaa_0 g)(x) \Eq g'(x) - x g(x),    \quad g\in \gg,
\] 
discussed above, whose properties are well documented: see, for example, 
Chen \& Shao~(2005, Lemmas~2.2 and 2.3).  Rather than take
the standard normal for~$\pi_0$, we actually prefer to perturb from a normal 
distribution $\nn(0,(1-\psi)^{-1})$.  This has
Stein operator 
\eq\label{Stein-normal}
  (\aaa_0 g)(x) \Eq g'(x) - (1-\psi)xg(x), \quad g\in \gg,
\en
which
% , choosing \[  w\ :=\ \{1-\half\psi(1-\tfrac1m)\}, \]   
gives
\[
  (\uu g)(x) \Eq -x\,
% \Blb 1-\psi-w + 
\frac{\psi(m+1)}{m+x^2}\, g(x),\quad  g\in \gg.
\]
The properties of~$\aaomi$ are as given in \prax,
with~$y$ replaced by~$g$.  For the supremum norm on~$\ff$, we find that
assumptions \Ref{norm-cond1}--\Ref{Stein-operator0-3} and 
\Ref{Stein-operator1-1}--\Ref{Stein-operator1-2} are satisfied, and that
\eqs
   \sup_x |x(\aaomi P_0 f)(x)| &\le& 2(1-\psi)^{-1}\|f\|;  \\
   \|\uu\aaomi P_0\| &\le& 2 \ps(1-\psi)^{-1}(1 + \tfrac1m)\ =:\ \g,
\ens
from \prax\,(b)(iii).
\adb{Condition~\Ref{key} is satisfied if $\g<1$, in which case 
Theorem~\ref{very-useful} shows that Stein's method works.}

Note, however, that 
Student's $t_m$~distribution itself is too far from the normal for this perturbation
argument to be applied, since then $\psi=1$, and so $\g=\infty$.

% As observed at the end of the previous section, 
% The choice of the set~$\hh$
% of test functions is critical, when working with the supremum norm. However,
For bounded functions~$f$ with bounded derivative, it 
follows from \prax\,(c)(iv) that
\[
  \sup_x |x(\aaomi P_0 f)'(x)| 
     \Le \frac3{1-\psi}\,\|f'\nsup\,.
\]
This translates into a bound for $\|\uu\aaomi P_0 f\|\ui$,
% Hence $\uu\aaomi P_0$ maps $\bff$ to~$\bff$, 
% if we now take 
% similar considerations could also have been applied, taking 
% the norm  $ \|f\| \ :=\ \|f\nsup + \|f'\nsup$ on~$\ff$, as in \Ref{diff-norm},
and~\Ref{key} is then satisfied for all~$\psi$ small enough. 
As for normal approximation, bounding $\ex\{(\aaa_0 g)(W)\}$ by a
linear combination of $\|g\nsup$,  $\|g'\nsup$ and~$\|g''\nsup$
may be a much more reasonable prospect than using only
$\|g\nsup$ and  $\|g'\nsup$, and
these quantities are themselves all bounded by multiples of
$\|f\|\ui$,
% as defined in~\Ref{diff-norm}, 
for $g=\aaomi P_0 f$ and $f\in\bff\ui:=\{f\in\ff\colon \|f\|\ui < \infty\}$:
see \prax\,(c)(i)--(iii), with $y=g$. In such cases, $d\ui$-approximation
is a consequence. 

To deduce Kolmogorov distance using Theorem~\ref{K-distance}, 
note that, for $\hh = \hh^K$, 
\eqs
  \g_{\hh} &=& \sup_{h\in \hh^K}\|\uu g_h^0\|\ui \\
	&\le& 2 \frac{\psi}{1-\ps}\Bl1 + \frac1m\Br
	   \Blb 1 + \frac1{\sqrt m} + \frac14 \,\sqrt{2\pi(1-\ps)}
		   + \frac12(1-\ps)\sqrt m \Brb,
\ens
from \prax\,(a)(i)--(iii).
If an approximation with respect to~$d\ui$ can be obtained from 
Theorem~\ref{very-useful}, then the estimate used in~\Ref{Stein-approx}
can be used also in~\Ref{cond-2}, and
the main further obstacle is thus to verify condition~\Ref{cond-1}.

\bigskip 
\nin{\bf Example~\ref{examples}.2}.\ 
Our second example also concerns a perturbation of the normal distribution, but
now to a distribution~$\pi_1$, whose Stein operator is not so easy to handle
directly.  This time, we take for~$\gg$ the space of real functions~$g$ with
$g(0)=0$ and having
bounded first and second derivatives, endowed with the norm
\[
  \|g\ngg \ :=\ \|g'\nsup + \|g''\nsup.
\]
As Stein operator~$\aaa_1$, we fix $\a>0$ and take the expression
\eq\label{2-St-op}
  (\aaa_1 g)(x) \Eq g''(x)  - xg'(x) + \a\{g(x+z) - g(x)\}, \quad g\in\gg,
\en
which can be viewed as a perturbation of the Stein operator
\[
   (\aaa_0 g)(x) \Eq g''(x) - xg'(x), \quad g\in\gg,
\]
characterizing the standard normal distribution. This operator is
equivalent to that given in~\Ref{Stein-standard}, and the properties
of~$\aaomi$ are given in \prax, with $y=g'$ and $\ps=0$.
The distribution~$\pi_1$ is
that of the equilibrium of a jump--diffusion process~$X$, with unit
infinitesimal variance, and having jumps of size~$z$ at rate~$\a$.

Once again, taking the supremum norm on~$\ff$, it is easy to
check that assumptions \Ref{norm-cond1}--\Ref{Stein-operator0-3} and 
\Ref{Stein-operator1-1}--\Ref{Stein-operator1-2} are satisfied, and since
\[
   \|(\aaomi P_0 f)'\nsup \Le \derivbnd\, \|f\nsup,
\]
from \prax\,(b)(i), it follows that 
\eq\label{2-gamma-sup}
    \|\uu\aaomi P_0\| \Le \derivbnd\,z\a.
\en
\adb{Hence, from Theorem~\ref{very-useful}, Stein's method works for~$\pi_1$ if
$\g = \derivbnd z\a < 1$;   
an estimate of the form~\Ref{Stein-approx} is all that is needed.}

As above, the supremum norm may be more
difficult to exploit in practice than the norm~$\|\cdot\|\ui$.
Here, for $f\in\bff$, and writing $g_f^0 = \aaomi P_0 f$, we have
\[
   |(\uu g_f^0)'(x)| \Le \a\int_0^z |(g_f^0)''(x+t)|\,dt \Le 4\a z \|f\nsup,
\]
from \prax\,(b)(ii), and Theorem~\ref{very-useful} can be applied if~$\a$
is small enough that $\g = (4+\derivbnd) z\a < 1$. 

For Kolmogorov approximation, note that, for $\hh = \hh^K$, 
\[
  \g_{\hh} \Eq \sup_{h\in \hh^K}\|\uu g_h^0\|\ui \Le (1+\derivbnd/4)\,z\a\,,
\]
by \prax\,(a)(i)--(ii).
Once again, the main effort in addition to~$d\ui$--approximation is to
verify~\Ref{cond-1} of Theorem~\ref{K-distance}.

Note that we are also free to perturb from other normal distributions.
If we choose to centre at the mean~$\a z$ of~$\pi_1$, we can do so
by writing
\[
  (\aaa_1  g)(x) \Eq g''(x)  - (x-\a z)g'(x) + \a\{g(x+z) - g(x) - zg'(x)\}, 
     \quad g\in\gg,
\]
with the first two terms the Stein operator for the normal distribution
$\nn(\a z,1)$.  The third, perturbation term
can be bounded by $2\a z^2 \|f\nsup$, and its derivative by
 $\a z^2\|f'\nsup$ (\prax\,(b)(ii)--(iii)), enabling
\Ref{key} to be satisfied for $\|f\|\ui$ for
 a larger range of~$\a$, if~$z$ is small
enough.  It is also possible to begin with $\nn(\a z, 1 + \a z^2/2)$,
correcting for both mean and variance.

It is also possible to generalize the class of perturbed measures by
replacing the term $\a(g(x+z)-g(x))$ corresponding to Poisson jumps
of rate~$\a$ and magnitude~$z$ by a more general L\'evy process,
taking instead $\int\{g(x+z)-g(x)\}\,\a(dz)$, for a suitable measure~$\a$.

\bigskip 
\nin{\bf Example~\ref{examples}.3}.\ 
As our third example, considered already in \BX\ (1999) and in \BC~(2002),
we consider (signed) compound Poisson distributions~$\pi_1$ on~$\integ$, 
the set of all integers, as 
perturbations of Poisson distributions on~$\integ_+:=\{0,1,2,\cdots\}$.  
We begin with~$\pi_1$
as the compound Poisson distribution CP$(\l,\mu)$ on~$\integ_+$, the
distribution of $\sli lN_l$, where $N_1,N_2,\ldots$ are independent,
and $N_l \sim \Po(\l \mu_l)$; $m_1 := \sli l\mu_l$ is assumed to be finite. 
In this case, we have $\xx = \xx_0 = \integ_+$.
With~$\gg$ the space of bounded functions $g\colon\nat\to\re$,
endowed with the supremum norm, a
suitable Stein operator for~$\pi_1$ is given by
\eq\label{CP-op}
  (\aaa_1 g)(j) \Eq \l\sli l\mu_l g(j+l) - jg(j), \quad j\ge0,
\en
considered as a perturbation of the Stein operator
\eq\label{Po-op}
  (\aaa_0 g)(j) \Eq \l m_1 g(j+1) - jg(j), \quad j\ge0;
\en
this means that
\eq\label{CP-U}
  (\uu g)(j) \Eq \l\sli l\mu_l\{g(j+l) - g(j+1)\}, \quad j\ge0.
\en     
Taking the supremum norm on~$\ff$, assumptions 
\Ref{norm-cond1}--\Ref{Stein-operator0-3} and 
\Ref{Stein-operator1-1}--\Ref{Stein-operator1-2} are satisfied; and since,
from the well-known properties of the solution of the Stein Poisson
equation,
\eq\label{ADB-alpha}
   \|\D(\aaomi P_0f)\nsup \Le \frac2{\l m_1}\,\|f\nsup,
\en
where $\D g(j) := g(j+1)-g(j)$, it follows that
\[
   \|\uu\aaomi P_0\| \Le 2\l\sli l(l-1)\mu_l /(\l m_1) \Eq 2m_2/m_1,
\]
where $m_2 = \sli l(l-1)\mu_l$.  Hence~\Ref{key} is satisfied if
$m_2/m_1 < 1/2$, and Theorem~\ref{very-useful} 
can then be invoked. 
Note that, in this setting, it is reasonable to work in terms of the
supremum norm, since total variation approximation may genuinely
be accurate.

There are nonetheless other distances that are useful.  Two such are
the Wasserstein distance~$d_W$, defined for  measures $P$ and~$Q$
on $\integ$ by
\[
  d_W(P,Q) \ :=\ \sup_{f\in\Lip_1} |P(f) - Q(f)|,
\]
where $\Lip_1 := \{f\colon \integ \to \re;\,\|\D f\nsup \le 1\}$, and the point
metric~$\dpt$ defined by
\[
  \dpt(P,Q) \ :=\ \max_{j\in\integ} |P\{j\} - Q\{j\}|,
\]
which has application when proving local limit theorems.

For Wasserstein distance, it is natural to begin with the semi-norm
$\|f\| := \|f\|_W := \|\D f\nsup$ on~$\ff$, which becomes a norm when restricted
to~$\bff'$.  The arguments in Section~\ref{formal} go through in this
modified setting very much as before; the only practical differences
are that one needs to check that $P_0\uu\aaomi$ maps $\bff'_0$ into
itself, and to replace the condition~\Ref{key} by
\eq\label{key'}
  \g \ :=\ \|P_0\uu\aaomi\| \ <\ 1.
\en
For the Poisson operator~$\aaa_0$  given in~\Ref{Po-op},  it is known
that
\eqa
   &&\|g_f^0\nsup \Le \|P_0f\|_W \Eq \|f\|_W;\quad 
	 \|\D g_f^0\nsup \Le 1.15(\l m_1)^{-1/2}\|f\|_W;\phantom{HHH}\non\\
   &&\|\D^2 g_f^0\nsup \Le 2(\l m_1)^{-1}\|f\|_W, \label{ADB-beta}
\ena
whenever $f\in\bff$ and $g_f^0 := \aaomi P_0 f$ 
[Barbour and Xia~(2005)].  Hence,
for~$\aaa_1$ as in~\Ref{CP-op} and $f \in \bff'_0$, it follows from~\Ref{CP-U} 
that
\[
  \|P_0\uu g_f^0\|_W \Eq \|\D P_0\uu g_f^0\nsup \Le \l\sli l(l-1)\mu_l \|\D^2 g_f^0\nsup
  \Le 2 (m_2/m_1) \|f\|_W,
\]
so that $P_0\uu\aaomi$ indeed maps $\bff'_0$ into itself, and
$\g = \|P_0\uu\aaomi\| \le 2m_2/m_1$.  Thus~\Ref{key'} is
satisfied for $m_2/m_1 < 1/2$, and the perturbation approach
can then be invoked.

For the point metric, we take the $l_1$--norm $\|f\| := \|f\|_1 := \sjz|f(j)|$ 
on~$\ff$.
For $f\in\bff_0$ and $g_f^0 = \aaomi P_0 f$, we have
\eq\label{ADB-gamma}
  \|g_f^0\nsup \Le (\l m_1)^{-1}\|f\|_1;\qquad 
	\|\D g_f^0\|_1 \Eq \sji|\D g_f^0(j)| \Le 2(\l m_1)^{-1}\|f\|_1;
\en
both inequalities are consequences of the proof of the second inequality
in \BHJ~(1992, Lemma~1.1.1).
Hence, from~\Ref{CP-U}, it follows immediately that
\eqs
  \|\uu g_f^0\|_1 &=& \sjo|(\uu g_f^0)(j)| \\
  &\le& \l \sli l\mu_l \sjo \sum_{s=1}^{l-1} |\D g_f^0(j+s)|\\
  &\le& 2\l m_2(\l m_1)^{-1}\|f\|_1 \Eq 2(m_2/m_1)\|f\|_1,
\ens
so that condition~\Ref{key} is once again satisfied if $m_2/m_1 < 1/2$.

If, more generally, $\pi_1$ is a (signed) compound measure on~$\integ$,
with characteristic function 
\[
   \exp\Blb \l\sum_{l\in\integ}\m_l(e^{il\th} - 1) \Brb,
\]
similar considerations can be applied. Here, we now have
$\xx = \integ$, but~$\xx_0$ is still~$\integ_+$. The corresponding Stein
operator is formally exactly as in~\Ref{CP-op}, except that the $l$-sum now
runs over the whole of~$\integ$, and we require $m_1$ to be positive;
also, the role of~$m_2$ is now played by $m'_2 = \sum_{l\in\integ} l(l-1)|\mu_l|$. 
When applying Lemma~\ref{useful} and Theorem~\ref{very-useful}, 
we have the inequalities 
\[
   \k(\pi,\integ_-) \Le  2\,|\pi|(\integ_-)
\]
for use with~$d_{TV}$,
\[
  \k(\pi,\integ_-) \Le  \sum_{j<0} |\pi|\{j\}(|j| + \l)
\]
for~$d_W$, and, with the fact that $\max_j\pi_0(j)\le (2e\l)^{-1/2}$ 
[\BHJ~(1992, p.~262)],
\[
  \k(\pi,\integ_-) \Le  \frac1{\sqrt{2e\l}}\,|\pi|(\integ_-) + 
   \max_{l< 0}|\pi|\{l\} 
\]
for~$\dpt$.

\bigskip 
\nin{\bf Example~\ref{examples}.4}.\ 
In this example, the setting is similar to that in the preceding example,
but we now consider a compound Poisson distribution~$\pi_1 = {\rm CP\,}(\l^1,\mu^1)$ 
on~$\integ_+$ as a perturbation not of a Poisson distribution, but of another
compound Poisson distribution~$\pi_0 = {\rm CP\,}(\l^0,\mu^0)$ on~$\integ_+$.  
The reason for doing so is that  the solutions to the 
Stein equation are known to be well behaved only for rather restricted classes 
of compound Poisson distributions: 
see Barbour \& Utev~(1998), Barbour \& Xia~(2000).  The perturbation method
offers the possibility of expanding the class of those with good behaviour by
including neighbourhoods not only of the Poisson distributions, but also of
any other compound Poisson distributions whose Stein solutions can be controlled.
In particular, we shall suppose that the distribution~$\pi_0 = {\rm CP\,}(\l^0,\mu^0)$
is such that
\[
  j\mu_j^0 \Ge (j+1)\mu_{j+1}^0,\quad j\ge1,
\]
and that $\d := \mu_1^0 - 2\mu_2^0 > 0$, these conditions implying 
that, with $c_1(\l^0)=4-2(\d\l^0)^{-1/2}$ and 
$c_2(\l^0)=\thalf(\d\l^0)^{-1}+2\log^+(2(\d\l^0))$,
\eq\label{BCL}
  \|g_f^0\nsup \Le \{\d\l^0\}^{-1/2}c_1(\l^0) \|f\nsup \quad \mbox{and}\quad
  \|\D g_f^0\nsup \Le \{\d\l^0\}^{-1} c_2(\l^0)\|f\nsup,
\en
where, as usual, $g_f^0 := \aaomi P_0 f$; see \BCL~(1992, pp.~1854-5). Here, 
the Stein operators $\aaa_0$ 
and~$\aaa_1$ are given as in~\Ref{CP-op}, with the corresponding choices
of $\l$ and~$\mu$, giving
\[
  (\uu g)(j) \Eq \sli l\{\l^1\mu_l^1 - \l^0\mu_l^0\}g(j+l) , \quad j\ge0.
\] 
As in the previous example, we shall only consider
perturbations which preserve the mean, so that also
\[
  \l^1\sji j\mu_j^1 \Eq \l^0\sji j\mu_j^0.
\]

Taking the supremum norm on~$\ff$, assumptions 
\Ref{norm-cond1}--\Ref{Stein-operator0-3} and 
\Ref{Stein-operator1-1}--\Ref{Stein-operator1-2} are satisfied.
In order to express the contraction condition~\Ref{key}, write
\[
  E \ :=\ \tfrac12 \sli l |\l^1\mu_l^1 - \l^0\mu_l^0|,
\]
and define probability measures $\r$ and~$\s$ on~$\nat$ by
\[
  \r_l \Eq E^{-1} l(\l^1\mu_l^1 - \l^0\mu_l^0)^{+};\quad
  \s_l \Eq E^{-1} l(\l^0\mu_l^0 - \l^1\mu_l^1)^{+},\qquad l\ge1;
\]
set $\th := E d_W(\r,\s)$, where~$d_W$ denotes the Wasserstein
distance.  Then, using~\Ref{BCL}, it follows easily that
\[
  \|\uu\aaomi P_0\nsup \Le  \{\d\l^0\}^{-1} c_2(\l^0)\th \ =:\ \g,
\]
with~\Ref{key} satisfied if $\g < 1$.

\bigskip 
\nin{\bf Example~\ref{examples}.5}.\ 
In our last example, we consider solving the Stein equation for 
a point process, whose distribution~$\pi_1$ is close to that of a 
spatial Poisson process.  Let~$\xxo$ be a compact metric space,
% with metric~$d_0$ bounded by~$1$,
and let $\xx$ denote the
space of Radon measures (point configurations) on~$\xxo$.
Then a Poisson process on~$\xxo$ with intensity measure~$\L$
satisfying $\l := \L(\xxo) < \infty$ is a random element 
$\sum_{l=1}^N \d_{X_l}$ of~$\xxo$, where $N,X_1,X_2,\ldots$ are
all independent, $N \sim \Po(\l)$ and $X_l \sim \l^{-1}\L$ for $l\ge1$,
and~$\d_x$ denotes the unit mass at~$x$.
Its distribution~$\pi_0$ can be characterized by the fact that
$\pi_0(\aaa_0 g) = 0$ for all $g$ in 
\[
  \gg := \{g\colon \xx\to \re;\, g(\emptyset)=0,\,\|\D g\nsup < \infty\}, 
\]
where
\[
  (\aaa_0 g)(\xi) \ :=\ \int_\xxo\Blb (g(\xi+\d_x)-g(\xi))\L(dx) +
   (g(\xi-\d_x)-g(\xi))\xi(dx) \Brb,
\]
and $\|\D g\nsup := \sup_{\xi\in\xx, x \in \xxo}|g(\xi+\d_x) - g(\xi)|$.
Note that the Stein operator~$\aaa_0$ is the generator of a spatial 
immigration--death
process, with~$\pi_0$ as its equilibrium distribution.  For the measure~$\pi_1$,
we take the equilibrium distribution of another spatial immigration--death
process on~$\xxo$, with generator~$\aaa_1$ given by
\[
  (\aaa_1 g)(\xi) \ :=\ \int_\xxo\Blb (g(\xi+\d_x)-g(\xi))\L_1(\xi,dx) +
   (g(\xi-\d_x)-g(\xi))\xi(dx) \Brb;
\]
here, the immigration measure is allowed to depend on the current 
configuration~$\xi$.  We can write $\aaa_1 = \aaa_0 + \uu$ if we set
\[
  (\uu g)(\xi) \ :=\  \int_\xxo (g(\xi+\d_x)-g(\xi))(\L_1(\xi,dx) - \L(dx)) ,
\]
and we note that  $\xx_0 := \supp(\pi_0) = \xx$.

We begin by considering perturbations appropriate for total variation 
approximation, taking the set of functions $\ff\colon \xx \to \re$ with
the supremum norm~$\|\cdot\nsup$. Then, as in Barbour 
\& Brown~(1992, pp.~~12--13),
it is possible to define a right inverse~$\aaomi$ satisfying \Ref{Stein-operator0-2}
and~\Ref{Stein-operator0-3}, with $A = 2$.  To check that $\uu\aaomi P_0\colon
\bff\to\bff$, 
% so that we need only to check the condition~\Ref{key} in order to
% be able to apply Theorem~\ref{useful}.  Now, 
we combine the definition of~$\uu$ and~\Ref{Stein-operator0-3} to give
\[
  |(\uu g_f^0)(\xi)| \Le 2\LL_\xi(\xxo)\|f\nsup,
\]
where $\LL_\xi(\cdot)$ denotes the absolute difference between the measures
$\L_1(\xi,.)$ and~$\L(\cdot)$; hence we shall need in addition to assume
that
\eq\label{PP-cond}
  \tl \ :=\ 2\sup_{\xi\in\xx} \LL_\xi(\xxo) \  <\ \infty,
\en
in order to make progress.  If we do, then~\Ref{key} is satisfied with
$\g = \tl$ if $\tl < 1$, and \adb{ Theorem~\ref{very-useful} can be used
to show that Stein's method works.}

Total variation is often too strong a metric for comparing point process
distributions, and so an alternative metric~$d_2$ is proposed in Barbour 
\& Brown~(1992), based on test functions Lipschitz with respect to a metric
on~$\xxo$ which is bounded by~$1$. Similar
calculations can be carried out in this setting also; the condition needed
to satisfy~\Ref{key} is somewhat more stringent.  

Even the contraction condition~$\tl < 1$ is rather restrictive.  Consider
a hard-core model, in which $\xxo \subset \re^d$ has 
volume~$\vartheta$, 
$\L(dx) = dx$, and $\L(\xi,dx) = I[\xi(B(x,\e)) = 0]\,dx$, \adb{where~$I[C]$
denotes the indicator of the event~$C$; this specification of~$\L(\xi,dx)$ is such
that no} immigration is allowed within distance~$\e$ of a point of the
current configuration~$\xi$.  Then $\tl = \vartheta$, and contraction is only
achieved if the expected number~$\vartheta$ of points under~$\pi_0$ 
is less than~$1$.
However, one could also consider a model with~$\pi_1$ the equilibrium
distribution of a slightly different immigration death process, in which
\[
  \L(\xi,dx) \Eq \max\{I[\xi(B(x,\e)) = 0],I[\xi(\xxo) > 2\vartheta]\}dx;
\]
for large~$\vartheta$, the difference between the equilibrium distributions of
the two processes is small, but, for the new process, $\tl \le 2\vartheta a(\e)$,
where $a(\e)$ is the area of the $\e$-ball, meaning that models with
much larger expected numbers of points can still satisfy the contraction
condition. Nonetheless, these are still models in which, at distance~$\e$,
little interaction can be expected; the mean number of pairs of points
closer than~$\e$ to one another under~$\pi_0$ is about $\tfrac12 \vartheta a(\e)$,
and, if the contraction condition is satisfied, this has to
be less than~$\tfrac14$.

\section{Illustrations}\label{illustration} 
\setcounter{equation}{0} 
In this section, we illustrate how the perturbations described above can
be used in specific examples.

\bigskip
\nin{\bf Example~\ref{illustration}.1}.\ 
In our first illustration, we return to the setting and notation 
of Example~\ref{examples}.2,
and consider approximation by the distribution~$\pi_1$ whose Stein operator 
is given in~\Ref{2-St-op} above.  The distribution we wish to approximate is 
the equilibrium distribution~$\pi$ of another jump--diffusion process, in which 
the jumps do not have fixed size~$z$, but are randomly chosen with~$z$ as mean;
this process has generator~$\aaa$ given by
\eq\label{2-St-op'}
  (\aaa g)(x) \Eq g''(x)  - xg'(x) + \a\int\{g(x+\z) - g(x)\}\mu(d\z), \quad g\in\gg.
\en
This distribution can be expected to be close to~$\pi_1$ provided
that the probability distribution~$\mu$ is concentrated about~$z$,
and since the distribution~$\pi$ is reasonably well understood, such an
approximation may constitute a useful simplification.

The main step is thus to establish a bound of the form~\Ref{Stein-approx}, after
which Theorem~\ref{very-useful} can be applied.  However, for $X\sim\pi$
and $g_f^0 := \aaomi P_0 f$, we immediately have
\eqs
  \pi(\aaa_1 g_f^0) &=& \pi(\aaa g_f^0) - \a \ex\Blb \int\{g_f^0(X+\z) - g_f^0(X+z)\}\mu(d\z)\Brb\\
   &=&  - \a \ex\Blb \int\{g_f^0(X+\z) - g_f^0(X+z)\}\mu(d\z)\Brb,
\ens
since $\pi(\aaa g) = 0$ for all $g\in\gg$. From this it follows by the mean value
theorem that
\[
  |\pi(\aaa_1 g_f^0)| \Le \half\a \int (\z-z)^2 \mu(d\z)\,\|(g_f^0)''\nsup 
     \Le 2\a \int (\z-z)^2 \mu(d\z)\,\|f\nsup.
\]
This suggests that the supremum norm on~$\ff$ is an appropriate choice, and
from Theorem~\ref{very-useful}, if $\g := \derivbnd\,z\a < 1$,
as in~\Ref{2-gamma-sup}, it follows that 
\[
   d_{TV}(\pi,\pi_1) \Le 2\a \int (\z-z)^2 \mu(d\z)/(1-\g).
\]
Thus the total variation distance between the two distributions is small if the
variance of~$\mu$ is small (and $\g < 1$).

\bigskip
\nin{\bf Example~\ref{illustration}.2}.\ 
We continue with the setting and notation of Example~\ref{examples}.2,
and again approximate by the distribution~$\pi_1$.  Here, as the measure~$\pi$, 
we take the
equilibrium distribution of a Markov jump process~$W_N$, defined as follows.
We let $X_N$ be the pure jump Markov process on $\integ_+$
with transition rates given by
\eqs
   j &\to& j+1\quad\mbox{at rate}\quad N;\qquad j \ \to\ j-1\quad\mbox{at rate}\quad j;\\
   j &\to& j + \lfloor z\sqrt N\rfloor \quad\mbox{at rate}\quad \a,
\ens
and we then set $W_N(t) := \{X_N(t) - N\}/\sqrt N$. If~$z=0$, the equilibrium distribution 
of~$X_N$ is the Poisson distribution with mean~$N$, and that of~$W_N$ the centred and
normalized Poisson distribution, which is itself, for large~$N$, close to the normal 
in Kolmogorov distance, but not in total variation.  Here, we wish to find bounds 
for the accuracy of approximation by~$\pi_1$ when $z>0$.  As above, we need
a bound of the form~\Ref{Stein-approx}, so as to be able to apply 
Theorem~\ref{very-useful}.

Much as above, we begin by observing that $\pi(\aaa g) = 0$ for all $g\in\gg$, where
now, writing $w_{jN} := (j-N)/\sqrt N$ and $\h_N := 1/\sqrt N$, we have
\eqs
  &&(\aaa g)(w_{jN}) \Eq N\{g(w_{jN}+\h_N) - g(w_{jN})\}\\
  &&\qquad\mbox{}  + j\{g(w_{jN}-\h_N) - g(w_{jN})\} 
% \\ &&\qquad\quad\mbox{}    
    + \a\{g(w_{jN}+\lfloor z\sqrt N\rfloor\h_N) - g(w_{jN})\}.
\ens
Subtracting $(\aaa_1 g)(w_{jN})$ and using Taylor's expansion, 
it follows that
\eqs
  \lefteqn{|(\aaa g)(w_{jN}) - (\aaa_1 g)(w_{jN})|}\\
  &&\Le N^{-1/2}(\tfrac13\|g'''\nsup + \tfrac12 |w_{jN}|\|g''\nsup + \a\|g'\nsup),
\ens   
so that
\eq\label{4.2-estimate}
  |\pi(\aaa_1 g)| \Le 
  N^{-1/2}(\tfrac13\|g'''\nsup + \tfrac12 \ex|W_N|\|g''\nsup + \a\|g'\nsup).
\en
Note that, taking $g(w)=w$ and $g(w) = w^2$ respectively in $\pi(\aaa g)=0$, 
as we may, by Hamza \&
Klebaner~(1995, Theorem~2), it follows that $|\ex W_N|\le\a z$ and
\[
  2\ex\{W_N^2\} \Le 2N\h_N^2 + \h_N|\ex W_N| + 2\a z|\ex W_N| + \a z^2
	  \Le 2+\a z\h_N+2\a^2z^2+\a z^2,
\]
which implies that
\[
  \{\ex|W_N|\}^2 \Le  \ex\{W_N^2\} \Le 
  1+\thalf\a z\h_N+\a^2z^2+\thalf\a z^2;
\]
thus $\ex|W_N|$ is uniformly bounded in~$N$.
Furthermore, for $g = g_f^0 := \aaomi P_0 f$ and $f\in\bff\ui$, 
we can control the first three derivatives of~$g_f^0$ by using
\prax\ with $y=(g_f^0)'$, so that~\Ref{4.2-estimate} yields a bound of the form
\[
  |\pi(\aaa_1 g_f^0)| \Le CN^{-1/2}\|f\|\ui,
\]
for all $f\in\bff\ui$. In view of Theorem~\ref{very-useful}, this
translates into the bound
\[
    d\ui(\pi,\pi_1) \Le CN^{-1/2}/(1-\g)
\]
if $\g < 1$, where now, for $\|\cdot\|\ui$, we have $\g = (4 + \derivbnd)z\a$,
as in Example~\ref{examples}.2.

If, instead, Kolmogorov distance is of interest, then the only obstacle is
to verify~\Ref{cond-1} of Theorem~\ref{K-distance}. For $g=g_h^0$, the
estimate given in~\Ref{4.2-estimate} is fine, except for the first term:
it is no longer possible to bound the difference
\[
  D_N(w) \ :=\ 
   N\{g(w+\h_N) - g(w) + g(w-\h_N)) - g(w)\} - g''(w)
\]
by $\tfrac13\h_N\|g'''\nsup$, since, for $h = h_a = \bone_{(-\infty,a]}$,
$g'''(a)$ is not defined.  However, it is clear that $|D_N(w)| \le 2\|g''\nsup$
for all~$w$, and that, for $|w-a| > \h_N$,
\[
  |D_N(w)| \Le \sup_{|x-w| \le \h_N} |g'''(x)|.
\]
Now, for $h=h_a$, taking $a>0$ without real loss of generality, we have
\eq\label{3rd-deriv-bnd}
  |g'''(x)| \Le C_1 + C_2 ae^{-a(a-x)}\bone_{(0,a)}(x),\quad x\neq a,
\en
for universal constants $C_1$ and~$C_2$, so that $g'''$ is well behaved except
just below~$a$.  The bound~\Ref{3rd-deriv-bnd} can then be combined with 
the concentration inequality
\[
  \pr[W_N \in [a,b]] \Le \{\tfrac12(b-a) + \h_N\}(\ex|W_N| + \a z),
\]
obtained by taking $g'' = \bone_{[a-\h_N,b+\h_N]}$ and 
$g'(w) = \int_{(b-a)/2}^w g''(t)\,dt$ for any $a\le b$ in $\pi(\aaa g)=0$, 
to deduce a bound $\ex|D_N(W_N)| \le CN^{-1/2}$, and hence Kolmogorov
approximation also at rate~$N^{-1/2}$.  Total variation approximation
is of course never good, since $\law(W_N)$ gives probability~$1$ to a
discrete lattice, and~$\pi_1$ is absolutely continuous with respect to
Lebesgue measure.

\bigskip
\nin{\bf Example~\ref{illustration}.3}.\ (Borovkov--Pfeifer approximation) 
Borovkov \& Pfeifer~(1996) suggested using a single $n$-independent 
infinite convolution of simple signed measures as a correction to the
Poisson approximation to the distribution 
of a sum of independent indicator random variables. Their approximation is
particularly effective in the case that they treated, the number of records
in~$n$ i.i.d.\ trials.  Here, the approximation is not as complicated as
it might seem, because the generating function of the correcting measure
can be conveniently expressed in terms of gamma functions.  
Its accuracy is then of order~$O(n^{-2})$, which is way
better than the $O(1/\log n)$ error in the standard Poisson approximation.
Their approach was extended to the multivariate case of
independent summands in \Ceka~(2002) and Roos~(2003).
Note also that Roos~(2003) obtained asymptotically sharp constants
in the univariate case.
In this example, by treating their approximating measure 
as a perturbation of the Poisson, as in
Example~\ref{examples}.3,  we investigate Borovkov--Pfeifer approximation to 
the distribution of the sum of {\it dependent\/} Bernoulli random variables.

Let $I_i$, $i\ge1$, be dependent Bernoulli~$\Be(p_i)$ random variables. 
Define $W=\sum_{i=1}^n I_i$, $W^{(i)}=W-I_i$, and let $\twi$ be a random 
variable having the conditional distribution 
of $W^{(i)}$ given $I_i=1$; that is, for all $k\in\integ_+$, 
$\pr(\twi =k)= \pr(W^{(i)}=k\giv I_i=1)$.
Let
$$
  \lambda=\sum_{i=1}^n p_i; \qquad 
  \eta_1=\sum_{i=1}^n \Blb \frac{p_i}{1-2p_i}\Brb \ex|\twi-W^{(i)}|.
$$
The Borovkov--Pfeifer approximation is defined to be the convolution
of the Poisson distribution~$\Po(\l)$ and the signed measure~$\BP$
determined by its generating function:
\eq\widehat {\BP}(z) \Eq \prod_{i=1}^\infty\Bigl\{
  \big(1+p_i(z-1)\big)\exp\left\{-p_i(z-1)\right\}\Bigr\}.
             \label{bp3}
\en
Using the fact that 
\eqa
 \lefteqn{e^{-p(z-1)}(1+p(z-1)) \Eq \exp\left\{\ln(1+pz/q)-\ln(1+p/q) - p(z-1)\right\}}\non\\
 &&\qquad=\ \exp\left\{\frac{p^2}q(z-1) + \sum_{l=2}^\infty \frac{(-1)^{l+1}}{l}
    \Bl\frac{p}{q}\Br^l(z^l-1) \right\},\phantom{HHHH}\label{xia-1}
\ena
where $q = 1-p$,
one can see that~$\BP$ is a signed compound Poisson measure, provided
that $\sn p_i^2 < \infty$. Note that $\sum_{i=1}^\infty p_i=\infty$ is allowed,
as is indeed the case for record values, when $p_i = 1/i$.

\begin{theorem}\label{BP}
Assume that $p_i<1/3$, $i\ge1$, that $\sii p_i^2<\infty$, and that
\eq
   \theta_1: \Eq \frac{m'_2}{m_1} \Eq \frac{\sn p_i^2(1-2p_i)^{-2}}{\l}\ <\ \frac{1}{2}.
  \label{bp5}
\en
Then
\eqa
 &&d_{TV}(\scrl(W),\Po(\l)*\BP) \Le \frac{2}{\l(1-2\theta_{1})}
    \bigg(\sum_{i=n+1}^\infty\frac{p_i^2}{(1-2p_i)^2}+\eta_1\bigg),\label{bp6}\\
 &&\dpt(\scrl(W),\Po(\l)*\BP) \non\\
 &&\quad\Le \frac{2}{\l(1-2\theta_{1})}
    \bigg(\sup_k\pr(W=k)\sum_{i=n+1}^\infty\frac{p_i^2}{(1-2p_i)^2}+\eta_1 \bigg),
    \hspace{1.1in}\label{bp7}\\
 &&d_W(\scrl(W),\Po(\l)*\BP) \Le \frac{1.15}{\sqrt{\l}(1-2\theta_{1})}
    \bigg(\sum_{i=n+1}^\infty\frac{p_i^2}{(1-2p_i)^2}+\eta_1\bigg).\label{bp8}
\ena
\end{theorem}

\nin{\bf Remark.} Let $I_i$, $i\ge1$, be independent. Then it suffices 
to prove the corresponding approximation for the sum
$W_s := \sum_{i=s}^nI_i$ only. Indeed, let~$\BP_s$ be specified by the 
generating function:
$$
  \widehat {\BP_s}(z) \Eq \prod_{i=s}^\infty
    \Bigl\{\big(1+p_i(z-1)\big)\exp\left\{-p_i(z-1)\right\}\Bigr\}.
$$
Then
\[
  \Po(\l)*\BP \Eq
    \law\Bl \sum_{i=1}^{s-1}I_i \Br * \Po\Bl\sum_{i=s}^n p_i\Br * \BP_s 
\]
and
\[
  \law(W) \Eq  \law\Bl \sum_{i=1}^{s-1}I_i \Br * \law(W_s),
\]
and so, by the properties of total variation we have
$$
  d_{TV}\left(\scrl\left(W\right),\Po(\l)*\BP\right)
     \le d_{TV}\left(\scrl\left(W_s\right),\Po\Bl\sum_{i=s}^n p_i\Br *\BP_s\right).
$$

If $W$ is the sum of independent Bernoulli variables, then $\h_1=0$ and
\[
    \sup_k \pr(W=k)\Le \Bl 4\,{\sn p_i(1-p_i)}\Br^{-1/2}, 
\]
see Barbour \& Jensen (1989, Lemma 1). Now, if we consider the records example 
of Borovkov \& Pfeifer~(1996), with $p_i=1/i$,
we can take any $s\ge4$ in the remark above, and obtain
orders of accuracy for the total variation distance, point metric and
Wasserstein metric of $O((n\ln n)^{-1})$, $O(n^{-1}(\ln
n)^{-3/2})$ and $O(n^{-1}(\ln n)^{-1/2})$, respectively. 

\bigskip
\nin{\bf Proof of Theorem~\ref{BP}.} In this case, $\xx=\xx_0=\integ_+$. 
Using \Ref{xia-1}, and setting $q_i = 1-p_i$, we can write~$\Po(\l)*\BP$ as 
the signed compound Poisson measure with generating function
\eq
  \exp\left\{\sli \l_l(z^l-1)\right\}, \label{bp14}
\en
where $\l_l=\l_{1l}+\l_{2l}$, with
\eqs
  \l_{1l}&=&\frac{(-1)^{l+1}}{l}\sum_{i=1}^n\left(\frac{p_i}{q_i}\right)^l, \quad l\ge1;\\
  \l_{21}&=&\sum_{i=n+1}^\infty\frac{p_i^2}{q_i}; \qquad
  \l_{2l}\Eq \frac{(-1)^{l+1}}{l}\sum_{i=n+1}^\infty\left(\frac{p_i}{q_i}\right)^l,
    \quad l\ge2.
\ens
Here, the components~$\l_{1l}$ come from the signed compound Poisson representation
of a sum of independent Bernoulli~$\Be(p_i)$ random variables, $1\le i\le n$,
and the~$\l_{2l}$ from the remaining~$\BP_{n+1}$ measure.

Let $\m_l=\l_l/\l$. Then, since 
$\sum_{l=1}^\infty l\l_{1l}=\sn p_i=\l$ and $\sum_{l=1}^\infty l\l_{2l}=0$,
we have $m_1=\sum_{l=1}^\infty l\m_l=1$. Hence,
the formula for~$\theta_1$ follows directly from
$$
  \sum_{l=2}^\infty l(l-1)|\l_l| 
    \Eq \sum_{l=2}^\infty(l-1)\sum_{i=1}^\infty\left(\frac{p_i}{q_i}\right)^l
    \Eq \sum_{i=1}^\infty p_i^2(1-2p_i)^{-2}.
$$

Next, we take Stein operators $\aaa_0$ as in \Ref{Po-op} and
$\aaa_1$ as in \Ref{CP-op}. For $g=g_f^0:=\aaa_0^{-1}P_0f$, it follows that
\eq
  \ex(\aaa_1 g)(W) \Eq  \Blb\sum_{l=1}^\infty l\l_{1l}\,\ex\,g(W+l)-\ex \{Wg(W)\}\Brb
    + \sum_{l=1}^\infty l\l_{2l}\,\ex\,g(W+l). \label{bp19}
\en
We begin by bounding the quantity in braces, which gives a bound for the accuracy
of the approximation of~$\law(W)$ by the distribution of a sum of independent
Bernoulli~$\Be(p_i)$ random variables.  We observe immediately that, for any~$i$
and~$l$,
\[
  \ex g(W+l) \Eq p_i \ex g(\twi+l+1) + q_i\ex\{g(W\uii+l) \giv I_i=0\}
\]
and that
\[
   \ex g(W\uii+l) \Eq p_i \ex g(\twi+l) + q_i\ex\{g(W\uii+l) \giv I_i=0\},
\] 
from which it follows that
\[
  \ex g(W+l) \Eq q_i \ex g(W\uii+l) + p_i \ex g(\twi+l+1) + p_i\uil,
\]
where we write $\uil := \ex g(W\uii+l) - \ex g(\twi+l)$.  Setting $\vil
:= (-1)^{l+1}(p_i/q_i)^l$, so that $l\l_{1l} = \sn \vil$, 
and observing that $p_i\vil = -q_i v_{i,l+1}$, we thus have
\eqs
  \lefteqn{\sli \vil \ex g(W+l) - \ex\{I_i g(W)\}}\\
  &&=\ q_i\sli\vil\ex g(W\uii+l) - q_i\sum_{l\ge2}\vil \ex g(\twi+l) \\
  &&\qquad\mbox{}  + p_i\sli \vil\uil - p_i\ex g(\twi+1)\\
	&&=\ \sli\vil\uil.
\ens
Adding over $1\le i\le n$, we thus find that
\eqs
  \lefteqn{\Blm \sum_{l=1}^\infty l\l_{1l}\,\ex g(W+l)-\ex \{Wg(W)\}\Brm}\\
   &&\Le \sn\sli |\vil|\,|\ex g(W\uii+l) - \ex g(\twi+l)| \Le \h_1\|\D g\|_\infty.
\ens      

\ignore{
By the identities that
\eqs
&&\D^0g(W+1)=g(W+1),\ g(W+l)=\sum_{k=1}^l\left(\begin{array}{c}l-1\\
k-1\end{array}\right)\D^{k-1}g(W+1),\\
&&\sum_{l=k}^\infty l\left(\begin{array}{c}l-1\\
k-1\end{array}\right)\l_{1l}=(-1)^{k+1}\sin p_i^k,
\ens
the first sum in (\ref{bp19}) is equal to
\eqs
&&\ex\sum_{l=1}^\infty
l\l_{1l}g(W+l)=\ex\sum_{l=1}^\infty
l\l_{1l}\sum_{k=1}^l\left(\begin{array}{c}l-1\\ k-1\end{array}\right)\D^{k-1}g(W+1)\\
&=&\ex\sum_{k=1}^\infty\sum_{l=k}^\infty
l\l_{1l}\left(\begin{array}{c}{l-1}\\ {k-1}\end{array}\right)\D^{k-1}g(W+1)=
\ex\sum_{k=1}^\infty\D^{k-1}g(W+1)(-1)^{k+1}\sin p_i^k\\
&=& \sum_{k=1}^\infty (-1)^{k+1}\sin p_i^k\Big(
p_i\ex\left\{\left.\D^{k-1}g(W^{(i)}+2)\right|I_i=1\right\}+q_i\ex\left\{\left.\D^{k-1}g(W^{(i)}+1)\right|I_i=0\right\}\Big)\\
&=& \sum_{k=1}^\infty(-1)^{k+1}\sin
p_i^k\Big(\ex\left\{\left.\D^{k-1}g(W^{(i)}+1)\right|I_i=1\right\}+p_i\ex\left\{\left.\D^{k}g(W^{(i)}+1)\right|I_i=1\right\}\Big)\\
&&+\sum_{k=1}^\infty(-1)^{k+1}\sin
p_i^kq_i\Big(\ex\left\{\left.\D^{k-1}g(W^{(i)}+1)\right|I_i=0\right\}-\ex\left\{\left.\D^{k-1}g(W^{(i)}+1)\right|I_i=1\right\}\Big)\\
&=& \sum_{k=1}^\infty(-1)^{k+1}\sin p_i^k\ex\left\{\left.\D^{k-1}g(W^{(i)}+1)\right|I_i=1\right\}\\
&&+\sum_{k=2}^\infty(-1)^k\sin p_i^k\ex\left\{\left.\D^{k-1}g(W^{(i)}+1)\right|I_i=1\right\}
\\
&&+\sum_{k=1}^\infty(-1)^{k+1}\sin
p_i^k\Big(\ex\D^{k-1}g(W^{(i)}+1)-\ex\left\{\left.\D^{k-1}g(W^{(i)}+1)\right|I_i=1\right\}\Big)\\
&=&\sin p_i\ex\big\{g(W^{(i)}+1)\vert I_i=1\big\}+
\sum_{k=1}^\infty(-1)^{k+1}\sin
p_i^k\Big(\ex\D^{k-1}g(W^{(i)}+1)-\ex\D^{k-1}g(\tilde
W^{(i)}+1)\Big).
\ens
Direct verification yields
$$\sin p_i\ex\big\{g(W^{(i)}+1)\vert I_i=1\big\}=\ex Wg(W)
$$
and
\eqs
&&\left|\sum_{k=1}^\infty(-1)^{k+1}\sin
p_i^k\bigg(\ex\D^{k-1}g(W^{(i)}+1)-\ex\D^{k-1}g(\tilde
W^{(i)}+1)\bigg)\right|\\
&\le& \sum_{k=1}^\infty\|\D^k g\|_\infty\sin
p_i^k\ex |W^{(i)}-\tilde W^{(i)}|\le \|\D
g\|_\infty\sin\ex|W^{(i)}-\tilde W^{(i)}|\sum_{k=1}^\infty
2^kp_i^k\\
&\le&\|\D g\|_\infty\eta_1.
\ens
}

It now remains to estimate the remaining element
$\sum_{l=1}^\infty l\l_{2l}\,\ex g(W+l)$ in~\Ref{bp19}. Using the identity 
$$
  g(W+l) \Eq g(W+1)+\sum_{s=1}^{l-1}\D g(W+s),
$$
we have
\eqs
  \sum_{l=1}^\infty l\l_{2l}\,\ex g(W+l)
    &=& \ex g(W+1)\Blb\sum_{l=1}^\infty l\l_{2l}\Brb
       + \sum_{l=1}^\infty l\l_{2l}\sum_{s=1}^{l-1}\ex\{\D g(W+s)\}\\
  &=& \sum_{l=1}^\infty l\l_{2l}\sum_{s=1}^{l-1}\ex\{\D g(W+s)\}, 
\ens
because $\sum_{l=1}^\infty l\l_{2l}=0$.  Now we have
\eqs
  &&\left|\ex\{\D g(W+s)\}\right|\ \leq\ \min\Blb\|\D g\|_\infty,
 %   \left|\ex\{\D g(W+s)\}\right|\ \leq\ 
    \,\|\D g\|_1\,\max_k \pr(W=k)\Brb,\\
  &&\sum_{l=2}^\infty l(l-1)|\l_{2l}| \Eq \sum_{i=n+1}^\infty\frac{p_i^2}{(1-2p_i)^2}.
  %\phantom{HHHHHHHHHHHHHHHHHHHHHHHHHHHHHHHHHHHHHHH}
\ens
The estimates \Ref{bp6}--\Ref{bp8} thus follow directly from the 
inequalities \Ref{ADB-alpha}, \Ref{ADB-beta} and~\Ref{ADB-gamma}
in Example~\ref{examples}.3.  \ep

 \section{Appendix}
 \setcounter{equation}{0}
Here, we collect various properties of the solution~$y$ to the equation
\eq\label{appx-1}
  y'(x) - (1-\ps)xy(x) \Eq h(x) - \bar h_\ps,  \quad x \in \re,
\en
for given~$h$ and $0\le \ps < 1$, where $\bar h_\ps = \ex h(N)$,
for $N \sim \nn(0,(1-\ps)^{-1})$.

\begin{proposition}\label{appx}
%\hfil\break
%\vskip-60pt
\nin
\begin{enumerate}
\item[{\rm (a)}] If $h = \bone_{(-\infty,z]}$ for any $z\in\re$, then
\eqs
 {\rm (i)}&& \|y\|_\infty \Le \frac14\sqrt{\frac{2\pi}{1-\ps}};\\
 {\rm (ii)}&&	 \|y'\|_\infty \Le 1;\\
 {\rm (iii)}&&  \sup_x |xy(x)| \Le \frac1{1-\ps}.
\ens
\item[{\rm (b)}] If $h$ is bounded, then
\eqs
 {\rm (i)}&& \|y\|_\infty \Le \sqrt{\frac{2\pi}{1-\ps}}\,\|h\|_\infty;\\
 {\rm (ii)}&&	 \|y'\|_\infty \Le 4\,\|h\|_\infty;\\
 {\rm (iii)}&&  \sup_x |xy(x)| \Le \frac2{1-\ps}\,\|h\|_\infty.
\ens 
\item[{\rm (c)}] If $h$ is uniformly Lipschitz, then
\eqs
 {\rm (i)}&& \|y\|_\infty \Le \frac{2}{1-\ps}\,\|h'\|_\infty;\\
 {\rm (ii)}&&	 \|y'\|_\infty \Le \frac4{\sqrt{1-\ps}}\,\|h'\|_\infty;\\
 {\rm (iii)}&&	 \|y''\|_\infty \Le \frac2{\sqrt{1-\ps}}\,\|h'\|_\infty;\\
 {\rm (iv)}&&  \sup_x |xy'(x)| \Le \frac3{1-\ps}\,\|h'\|_\infty.
\ens
\end{enumerate}
\end{proposition}

\begin{proof}
Equation~\Ref{appx-1} can be transformed, using the substitution
$x = w/\sqrt{1-\ps}$, into the equation with $\ps=0$ for the
standard normal distribution, for which the corresponding bounds
are mostly given in Chen \& Shao~(2005, Lemmas 2.2 and~2.3).
In particular, the bounds (a)(i)--(iii) follow directly from their
Equations (2.9), (2.8) and~(2.7), respectively; the bounds 
(b)(i)--(ii) from the proofs of their Equations (2.11) and~(2.12); 
and the bounds (c)(i)--(iii) from their Equations (2.11)--(2.13).

The bound (b)(iii) is easily deduced from the explicit expression
for the solution~$y$:  for instance, for $x>0$, we have
\[
  xy(x) \Eq - xe^{(1-\ps)x^2/2}\int_x^\infty e^{-(1-\ps)t^2/2}
	  (h(t) - \bar h_\ps)\,dt\,,
\]
immediately giving
\[
  |xy(x)| \Le x	\int_0^\infty e^{-(1-\ps)zx}\,|h(x+z) - \bar h_\ps|\,dz\,,
\]
from which (b)(iii) follows.

For (c)(iv), we argue only for $x<0$, since the proof for $x>0$
is entirely similar.  Noting that 
\[
  y''(x)-(1-\psi)xy'(x) \Eq (1-\psi)y(x)+h'(x)\,,
\]	
we obtain
\[
  y'(x) \Eq e^{\frac{(1-\psi)x^2}{2}}\int_{-\infty}^x
	    \{(1-\psi)y(t)+h'(t)\}\, e^{-\frac{(1-\psi)t^2}{2}}\,dt\,;
\]			
hence
\eqs
  |xy'(x)|&\le& \{(1-\psi)\|y\|_\infty+\|h'\|_\infty\}|x|e^{\frac{(1-\psi)x^2}{2}}
	    \int_{-\infty}^x e^{-\frac{(1-\psi)t^2}{2}}\,dt \\
  &\le& \|y\|_\infty+\|h'\|_\infty/(1-\psi)\,.
\ens
But now, from (c)(i) above, 	
$\|y\|_\infty \le \frac{2}{1-\ps}\,\|h'\|_\infty\,$.\ep
\end{proof}

\bigskip\bigskip
\nin {\bf Acknowledgement.}  We wish to thank a referee for extremely helpful
suggestions as to presentation.	 We also wish to acknowledge the generous support 
that has made this research possible; in particular from
Schweizerischer Nationalfonds Projekt Nr.\ 20-107935/1 and from
the ARC Centre of Excellence for Mathematics and Statistics of Complex Systems (ADB)
and from ARC Discovery Grant DP0209179 (VC).

\def\ac{{Academic Press}~}
\def\aap{{Adv. Appl. Prob.}~}
\def\ap{{Ann. Probab.}~}
\def\anap{{Ann. Appl. Probab.}~}
\def\jap{{J. Appl. Probab.}~}
\def\jws{{John Wiley $\&$ Sons}~}
\def\bny{{New York}~}
\def\ptrf{{Probab. Theory Related Fields}~}
\def\sp{{Springer}~}
\def\spa{{Stochastic Processes and their Applications}~}
\def\sv{{Springer-Verlag}~}
\def\tpa{{Theory Probab. Appl.}~}
\def\zw{{Z. Wahrsch. Verw. Gebiete}~}

\end{document}